%% file: trapzd8.tex
\documentclass[12pt]{article}
\scrollmode

\usepackage{color}

\usepackage{amsfonts}
\usepackage{amsmath}
\usepackage{amssymb}
\usepackage{graphicx}

\usepackage{mathrsfs}

\allowdisplaybreaks

\title{On the dynamics of trap models in $\Z^d$} 
\date{}
\author{L.~R.~G.~Fontes \footnote{IME-USP, Rua do Mat\~ao 1010, 05508-090
S\~ao Paulo SP,  Brazil, lrenato@ime.usp.br} \thanks{Partially
supported by CNPq grant 307156/2007-9, and FAPESP grant 2004/07276-2.}
\and
P.~Mathieu \footnote{CMI, 39 rue Joliot-Curie, 13013 Marseille, France,
pierre.mathieu@cmi.univ-mrs.fr}      }

\input{def54}

\begin{document}

\maketitle


\input{abs54}


\input{intro54}


\input{mod54}


\input{pre30}


\input{conv54}


\input{age54}


\input{str53}


\input{acks}




\input{bib30}
\end{document}

%% file: def54.tex
\def\eps{\varepsilon}

\def\P{{\mathbb P}}
\def\esp{{\mathbb E}}
\def\Z{{\mathbb Z}}
\def\N{{\mathbb N}}

\def\Q{{\mathbb Q}}

\def\1{{\mathbf 1}}

\def\a{\alpha}
  
\def\l{\lambda}

\def\x{{\cal X}}
\def\X{{\cal X}}
\def\Q{{\cal Q}}

\def\d{\delta}

\def\ttn{\tilde\tau^{(N)}}
\def\ups{\Upsilon}

\def\tk{\tau^{(k)}}
\def\tked{\tau^{(k,\eps,\d)}}
\def\ttned{\tilde\tau^{(N,\eps,\d)}}
\def\xk{X^{(k)}}
\def\xl{X^{(\ell)}}

\def\ckn{C^{(k,N)}}
\def\ck{C^{(k)}}
\def\ckned{C^{(k,N,\eps,\d)}}
\def\cked{C^{(k,\eps,\d)}}
\def\ikn{I^{(k,N)}}
\def\ik{I^{(k)}}
\def\ikned{I^{(k,N,\eps,\d)}}
\def\iked{I^{(k,\eps,\d)}}
\def\ykn{Y^{(k,N)}}
\def\yk{Y^{(k)}}
\def\ykned{Y^{(k,N,\eps,\d)}}
\def\yked{Y^{(k,\eps,\d)}}
\def\yoed{Y^{(1,\eps,\d)}}

\def\rkned{\Delta^{\!(k,N,\eps,\d)}}
\def\roned{\Delta^{\!(1,N,\eps,\d)}}

\def\rked{\Delta^{\!(k,\eps,\d)}}
\def\r1ed{\Delta^{\!(1,\eps,\d)}}
\def\ykne{Y^{(k,N,\eps)}}
\def\yke{Y^{(k,\eps)}}
\def\y1e{Y^{(1,N,\eps)}}
\def\ye{Y^{(\eps)}}

\def\tj{\hat T^{(j)}}

\def\bce{\tilde C^{(\eps)}}
\def\bvd{V^{(\delta)}}
\def\bved{\tilde V^{(\eps,\delta)}}
\def\bzd{Z^{(\delta)}}
\def\bzed{\tilde Z^{(\eps,\delta)}}
\def\vyed{\bar Y^{(\eps,\delta)}}
\def\bwd{W^{(\delta)}}
\def\bwed{\tilde W^{(\eps,\delta)}}
\def\bced{\tilde C^{(\eps,\delta)}}

\def\hye{\hat Y^{(\eps)}}

\def\bzd{Z^{(\d)}}
\def\ddye{\ddot Y^{(\eps)}}
\def\bryed{\breve Y^{(\eps,\delta)}}
\def\byed{\hat Y^{(\eps,\delta)}}

\def\bte{\tilde T^{(\eps)}}
\def\ze{\zeta^{(\eps)}}

\def\xe{x^{(\eps)}}
\def\yep{y^{(\eps)}}
\def\mue{\mu^{(\eps)}}

\def\t{{\mathfrak T}}
\def\teps{{\mathfrak T}^{(\eps)}}
\def\te{\tilde\eps}

\def\E{{\cal A}}

\def\B{{\cal B}}
\def\I{{\cal I}}

\def\RR{{\cal R}}
\def\rk{{\cal R}^{(k)}}
\def\rl{{\cal R}^{(\ell)}}

\def\tse{\tilde{\cal S}^{(\eps)}}
\def\hse{\hat{\cal S}^{(\eps)}}
\def\ZZ{{\cal Z}}

\def\x1{X_1^{(1)}}

\def\txt{\hat X}

\def\tx0{\tilde X^0}

\def\nn{\nonumber}

\def\={&=&}

\def\+{&+&}

\def\lf{\lfloor}
\def\rf{\rfloor}

\newtheorem{theo}{Theorem}
\newtheorem{prop}[theo]{Proposition}
\newtheorem{lm}[theo]{Lemma}
\newtheorem{cor}[theo]{Corollary}
\newtheorem{rmk}[theo]{Remark}

\def\beq{\begin{equation}}
\def\eeq{\end{equation}}
\newcommand{\bei}{\begin{itemize}}
\newcommand{\eei}{\end{itemize}}
\newcommand{\ben}{\begin{enumerate}}
\newcommand{\een}{\end{enumerate}}
\newcommand{\beqn}{\begin{eqnarray}}
\newcommand{\beqnn}{\begin{eqnarray*}}
\newcommand{\eeqn}{\end{eqnarray}}
\newcommand{\eeqnn}{\end{eqnarray*}}
\newcommand{\brm}{\begin{rmk}}
\newcommand{\erm}{\end{rmk}}

%% file: abs54.tex
\begin{abstract}

We consider trap models in $\Z^d$. These are stochastic processes in a random environment as follows.
The  environment is given by a family $\tau=(\tau_x\,,\, x\in\Z^d)$ of positive iid random variables 
in the basin of attraction of an $\a$-stable law, $0<\a<1$. 
Given $\tau$, our process is a continuous time Markov pure jump process, 
whose jump chain is a in principle generic random walk in $\Z^d$, $d\geq1$, independent of $\tau$, and
$\tau$ represents the holding time averages of the continuous time process. We may think of the sites of $\Z^d$ as traps, and
of $\tau_x$ as the depth of trap $x$. We are interested in the {\em trap process}, namely the process that associates 
to time $t$ the depth of the currently visited trap. Our first result is the convergence of the law of that process under suitable scaling. 
The limit process is given by the jumps of a certain $\a$-stable subordinator at the inverse of another $\a$-stable subordinator, correlated with
the first subordinator.
For that result, the requirements on the underlying random walk are $a)$ the validity of a law of large numbers for its range, and 
$b)$ the slow variation at infinity of the tail of the distribution of its time of return to the origin: 
they include all transient random walks as well as all random walks in $d\geq2$, and also many one dimensional random walks,
but not the simple symmetric case. 
We then derive {\em aging results} for our process, namely scaling limits for some two-time correlation functions of the process; 
a strong form of those results requires an assumption of transience, stronger than $a,\,b$ above. 
The scaling limit result mentioned above is an averaged result with respect to the environment. 
Under an additional condition on the size of the intersection of the ranges of two independent copies of the underlying 
random walk, roughly saying that it is small compared with size of the range, we derive a stronger scaling limit result, 
roughly stating that it holds in probability with respect to the environment. With that additional condition, we also strengthen 
the aging results, from the averaged version mentioned above, to convergence in probability with respect to the environment.

\end{abstract}


\smallskip


%% file: intro54.tex
\section{Introduction} 
\setcounter{equation}{0}

We begin with a more precise definition of random walks among random traps. These are constructed 
through the following two-step procedure. 
We first choose a probability measure on $\Z^d\setminus\{0\}$, say $\pi$, and let $\tau=(\tau_x\,,\, x\in\Z^d)$ be a collection of 
positive real numbers attached to the points of $\Z^d$. 
The {\it random walk in the trap environment $\tau$} 
is then the continuous time Markov process with values in $\Z^d$ that starts at the origin and 
 has generator 
\begin{equation}\label{eq:df0}
 {\cal L}^\tau f(x)=\frac 1{\tau_x}\sum_y (f(y)-f(x))\pi(y-x)\,,
\end{equation} 
so that 
when sitting at point $x\in\Z^d$, 
the process waits for an exponentially distributed time of mean $\tau_x$ and then jumps to point $x+y$ where $y$ 
is sampled from the distribution $\pi$. This procedure is then iterated with independent hopping times and jumps. 
$\tau$ will be taken at random so that the $(\tau_x\,,\, x\in\Z^d)$ is an independent identically distributed 
family of random variables whose common law belongs to the domain of attraction of a 
stable law of index $\alpha\in(0,1)$. The model thus defined is therefore an example of a random walk in a random 
environment. 
We denote by $\P$ the probability thus defined. More precisely, $\P$ is a probability measure on the product space 
$\Omega\times{\cal D}([0,\infty),\Z^d)$ where $\Omega=(0,\infty)^{\Z^d}$ is the space of {\it trap environments} 
and ${\cal D}([0,\infty),\Z^d)$ is the space of c\`adl\`ag trajectories from $[0,\infty)$ to $\Z^d$. The first marginal of 
$\P$ is of the form $Q^{\Z^d}$, where $Q$ belongs to the domain of attraction of a 
stable law of index $\alpha\in(0,1)$, and the conditional law of the second marginal given $\tau$ is the law of 
the random walk in the trap environment $\tau$. We use the notation $(\X_t\,,\, t\geq 0)$ for the canonical 
projections defined on ${\cal D}([0,\infty),\Z^d)$ that give the position of the random walker 
at times $(t\geq 0)$.

\begin {rmk} 
It is not difficult to see that if we choose $Q$ with compact support in $(0,+\infty)$ then the behaviour of the random
walk with traps is very similar to the random walk without traps. For instance, if $\pi$ is symmetric with finite support, one 
finds that the scaling limit of $\X$ under diffusive scaling is a Brownian motion. Fluctuations of the environment 
only affect the value of the effective diffusivity. In order to observe stronger slowing down effects, in particular aging, one has to 
choose heavy tailed $\tau$'s as we do here. 
\end{rmk}

This process is an example of a {\it trap model} in the spirit of J-P.~Bouchaud.  
 One important aspect of it is the lack of dependence of $\pi$ on $\tau$. (A class of models where there is such a dependence,
known as {\it asymmetric trap models}, have also been considered in the physics and mathematics literature. See below. Unless explicitely mentioned,
we do not discuss these models here.)
Such processes were initially introduced in the context of statistical mechanics as 
toy models for spin glasses and in order to illustrate the phenomenon of {\it aging}, see \cite{kn:Bou}, 
\cite{kn:BD} or \cite{kn:BCKM} for instance. 
In usual models of spin glasses, the Hamiltonian is a random Gaussian field of large variance. 
 At low temperature, it is natural to guess that the main contributions to the dynamics come from states of 
 low energy. 
As the statistical properties of extremes of 
 log Gaussian fields, the Gibbs factors in this context, are described by random variables 
with polynomial tail, the choice of a law in the basin of attraction of a stable law for $\tau$, 
 which plays a similar role in the simplified model,
is also natural. Note that the parameter $\alpha$ can then be interpreted as the
temperature, see \cite{kn:BBM}  and Subsection 3.2 in~\cite{kn:FL}. 

The aging property refers to the following phenomenon: as time increases, 
the process visits a larger and larger part of 
its state space and therefore increases its probability to find a location $x$ where $\tau_x$ is large. 
Since the time the process stays at location $x$ before jumping off is of order $\tau_x$, some slow down 
effect might take place. One way to measure how much the process is slow is to compute quantities of 
the form
 
\begin{equation}\label{eq:af}
\Pi(s,t)=\P(\X_r=\X_t\,,\, r\in[t,t+s])\,,
\end{equation}  
which are generally called {\it aging functions} in this context.
The Markov property implies that 
\begin{equation}\label{eq:af1}
\Pi(s,t)=\esp(e^{-s/{\tau_{\X_t}}})\,,
\end{equation}  
and thus we observe that a non trivial limit for $\Pi(s,t)$ as $s,t\to\infty$, with $s$ and $t$ related in a given way,
implies that, at large time $t$, 
$\tau_{\X_t}$ should be of order $s$, so that the (order of the) 'age' of the process can actually be approximately 
read from its position at large times. Thus, in order to describe aging, we are led to considering the asymptotics of 
the {\it age} (or  {\it trap}) {\it process} 
$\E=(\E_t=\tau_{\X_t}\,,\,t\geq 0)$. 

The first computations of J-P.~Bouchaud and D.~Dean in \cite{kn:Bou} and \cite{kn:BD} consisted in describing the asymptotics 
of trapped random walks on a large complete graph and in some appropriate scaling. 
Since then, 
the subject has developped into a rich mathematical theory.  
 Mathematical papers treating the model in the complete graph include~\cite{kn:BF} and~\cite{kn:FM}.
Although one motivation is 
certainly to understand the physicists' claims and prove aging for as realistic as possible 
models of spin glasses, see \cite{kn:ABG1},  \cite{kn:ABG2} and more recently \cite{kn:BBC}, 
it also turns out that trapping and aging effects 
also play a role in models without any connection to spin glass theory such as 
random walks with random conductances or random walks on Galton-Watson trees, see \cite{kn:BC},  \cite{kn:AFGH}.
The main strategy used in these papers has a strong potential theoretic flavour: for a given 
realisation of the trap environment $\tau$, one tries to identify, among the different points 
 $x$ with large $\tau_x$ which will be hit by the random walk. 
We refer to \cite{kn:AC2} for a presentation of this point of 
view in an abstract setting. One advantage of this approach is that it does not seem to require the state space 
to have many symmetries. It provides strong forms of aging properties that are valid for a given 
realisation of the traps. On the other hand this machinery is often quite heavy to use. 

As far as trap models on $\Z^d$ are concerned, excluding the asymmetric case (where $\pi$ depends on $\tau$
in a specific way, as mentioned above; see~\cite{kn:BC},~\cite{kn:C2},~\cite{kn:AC},~\cite{kn:M}), 
only the case of the simple symmetric random walk was investigated so far.
It corresponds to $\pi$ being the uniform law on the nearest neighbors of the origin. 
Then the paths of the process $\X$, i.e. the sequence of the different points visited by $\X$, is a symmetric 
nearest neighbor  random walk on $\Z^d$. 
The speed at which the process $\X$ moves i.e. the different hopping times at the successive locations 
are given by the environment $\tau$. 
The one-dimensional case happens to be special:
due to the strong recurrence properties of the simple symmetric random walk on $\Z$, the process localizes.  
This localization effect, aging properties and scaling limits are precisely 
described in \cite{kn:FIN}. The scaling limit is a singular diffusion now known under the name of {\it FIN}. 
In higher dimension $d\geq 2$, the scaling limit is known to be the so called 
{\it Fractional Kinetics} process: $d$-dimensional Brownian motion time changed by the inverse of an 
independent stable subordinator as proved in \cite{kn:C}, \cite{kn:ACM} and \cite{kn:AC1}, the strategy 
being similar to \cite{kn:AC2}. Besides a number of estimates on the Green kernel of simple symmetric random 
walk in $\Z^d$, the proof in the $d=2$ case involves rather sophisticated renormalization technics.
The same result also follows in $d\geq5$ as a particular case of results in~\cite{kn:M}, where a different approach is developed.

What do we do here? 
A first motivation of this paper is to derive aging properties for a more general class of random walks 
than the nearest neighbors case, in the form of an appropriate scaling limit of the age process, as suggested by our discussion above; 
see Theorem~\ref{teo:conv1} in Section~\ref{sec:mod} below. 
In doing so we hope to clarify which properties of the random walk are 
truly relevant for aging. Observe in particular that the usual recurrence versus transience dichotomy does not apply here 
 (as we can already conclude from the results of \cite{kn:AC1} for the simple symmetric case). 
As an outcome, we obtain a new proof of the scaling limit that applies to any genuinely $d$ dimensional 
random walk for $d\geq 2$. This proof is more conceptual than the approach previously used by other authors. 
Indeed we need to know very little about specific estimates for the transition probabilities or Green kernels. 
We also completely avoid the renormalization step, even in the $d=2$ case. It should also be mentioned that 
 Theorem~\ref{teo:conv1} is an annealed result with respect to the environment, as 
opposed to the analogue quenched result of~\cite{kn:AC1}. 
In particular it applies in situations where $\X$ does not have a quenched scaling limit,
see Remark~\ref{rmk:nec_cond} in Section~\ref{sec:str} below.  
It can however be strengthened  with little more effort, say, half way towards a quenched 
result, under a natural additional condition on the intersections of the ranges of independent copies of our process
(see Theorem~\ref{teo:str} in Section~\ref{sec:str} below).

Observe that in our general setting it does not make sense to look for scaling limits of the process $\X$ 
itself. Indeed the underlying random walk (with increments distributed according to $\pi$) may not have a 
non trivial scaling limit. 
 We choose then to focus on the age process $\E$, which is a natural object in the aging context, as we discuss next.
The expression of the aging function $\Pi(s,t)$, suggests that we should look at the limiting law 
of $\E_t/s$ to derive  an aging property. We actually provide a more complete answer by describing the 
scaling limit of the full process $\E$ thus establishing a fuller aging picture. 
This scaling limit is expressed as the value of the jump of 
some subordinator computed at the inverse of another subordinator. 
Interestingly, this scaling limit is universal, even if the scale on which the process $\E$ lives depends 
on the random walk (and in particular is linear if and only if the random walk is transient). 

The topology under which we are able to establish Theorem~\ref{teo:conv1} is
quite weak, though, due to the nature of the age process. Obtaining a scaling limit result for aging functions
like~(\ref{eq:af}) requires more work, done in Section~\ref{sec:age} 
(for~(\ref{eq:af}) and two other examples) under the stronger assumption of transience,
and in Subsection~\ref{ssec:stra} with the additional condition of Section~\ref{sec:str}; see
Theorems~\ref{teo:age} and~\ref{teo:ages}. Integrated forms of those results follow from Theorem~\ref{teo:conv1},
under the original assumptions, as discussed separately at the end of those sections.

The remainder of this paper is organized as follows. In Section~\ref{sec:mod} we have a detailed presentation of the model,
assumptions and one result (annealed scaling limit of $\E$), with some more discussion. In Section~\ref{sec:rw}, we discuss some preliminary results
on random walks that are used subsequently. In Section~\ref{sec:conv}, we prove the annealed result just mentioned, and in
Sections~\ref{sec:age} and~\ref{sec:str} we state and prove our further scaling limit results for some aging functions, and in a stronger than
annealed sense, as discussed above.

%% file: mod54.tex
\section{Model and results}
\label{sec:mod}
\setcounter{equation}{0}

As in the introduction, let $\pi$ be a probability measure on $\Z^d\setminus\{0\}$, and let 
$\tau=(\tau_x\,,\, x\in\Z^d)$ be a collection of 
positive numbers attached to the points of $\Z^d$ 
chosen as follows.
Let $Q$ be a probability measure on $(0,\infty)$ that belongs to the domain of attraction of a 
stable law of index $\alpha\in (0,1)$. In other words we assume that 
\begin{equation}\label{eq:tail}
 Q(u,\infty)=\ell(u)\,u^{-\alpha},\,u>0
\end{equation}
where $\ell$ is a slowly varying function at infinity. 
We choose for 
\begin{equation}\label{eq:tau}
\tau=\{\tau_x\,,\, x\in\Z^d\} 
\end{equation} 
a family of independent random variables 
with law $Q$. More precisely we endow the product space $\Omega=(0,\infty)^{\Z^d}$ 
with the law $\Q=Q^{\Z^d}$.

We consider the Markov generator 
\begin{equation}\label{eq:df}
 {\cal L}^\tau f(x)=\frac 1{\tau_x}\sum_y (f(y)-f(x))\pi(y-x)\,.
\end{equation} 
Let $P^\tau_x$ be the law 
of the Markov process $\X$ generated by ${\cal L}^\tau$ and started at $x$ 
on path space ${\cal D}([0,\infty),\Z^d)$. 
We recall that $(\X_t\,,\, t\geq 0)$ denotes the canonical 
projections on ${\cal D}([0,\infty),\Z^d)$. 
We define the age process (as in the previous section): 
\begin{equation}\label{eq:en}
 \E=(\E_t=\tau_{\X_t}\,,\,t\geq 0)\,.
\end{equation}

The so-called {\it annealed law} of the process $\X$ 
is the semidirect product measure on $\Omega\times{\cal D}([0,\infty),\Z^d)$ 
defined by 
\begin{equation}\label{eq:joint}
\P(A\times B)=\int_A d\Q(\tau) P^\tau_0(B),
\end{equation} 
where $A$ and $B$ are measurable subsets of $\Omega$ and ${\cal D}([0,\infty),\Z^d)$ 
respectively. 

In order to state our assumptions we introduce an auxiliary random walk:  
let $\xi_1,\xi_2,\ldots$ be iid $\Z^d$-valued 
random vectors with distribution $\pi$ and define 
\begin{equation}\label{eq:rw}
 X_0=0,\quad X_n=\sum_{i=1}^n\xi_i,\quad n\geq1\,.
\end{equation}
$X=(X_n,\,n\geq0)$ is a version of the jump chain of $\X$.

Let us also define the {\em range} of $X$ ({\em up to time} $n$), 
\begin{equation}\label{eq:range}
\RR_n=\RR_n(X)=\{z\in\Z^d:\,X_i=z\mbox{ for some }i\leq n\},\quad n\geq0,
\end{equation}
and make 
\begin{equation}\label{eq:range1}
R_n=|\RR_n| \mbox{ and } \rho_n=\esp(R_n). 
\end{equation}
Also define for $n\geq1$
\begin{equation}\label{rn}
r_n=\P(X_1\ne0,\ldots,X_n\ne0)\,.
\end{equation}
We will at times think of $(r_n)$ as a function.

Our first result requires the following assumptions.

\noindent{\bf Assumption} $\mathbf A$ (law of large numbers for the range): 
\begin{equation}\label{llnr}
 \lim_{n\to\infty} \frac {R_n}{\rho_n}=1\,\mbox{ in probability}. 
\end{equation} 
	    

\noindent{\bf Assumption} $\mathbf B$ (slow variation of $r$): $r:\N\to[0,1]$ given in~(\ref{rn}) above
is slowly varying at infinity. \label{assb}

\begin{rmk}
 \label{rmk:lln}
All transient random walks in $\Z^d$, $d\geq1$, including all random walks in $d\geq3$, 
obviously satisfy Assumptions $A$ and $B$.
But all planar random walks~\cite{kn:JP1, kn:JP2}, and $1$-dimensional  $\beta$-stable 
random walks with  $\beta\leq1$~\cite{kn:LR} also satisfy Assumptions $A$ and $B$.
\end{rmk}

Before stating our first result, we describe the form of the scaling limit of $\E$: we introduce an $\a$-stable subordinator 
$\ups=(\ups_r)_{r\geq0}$ and a family of independent mean 1 exponential random variables $\{T_r;\,r\geq0\}$, and let 
\begin{equation}
\label{eq:V}
V_s=\int_0^sT_r\,d\ups_r,\quad s\geq0,
\end{equation}
and $W=V^{-1}$ be the inverse of $V$ (see Remark~\ref{rmk:inv} below). Let finally
\begin{equation}
\label{eq:Z}
Z_t=\ups_{W_t}-\ups_{W_t-},\quad t\geq0.
\end{equation}

\begin{rmk}
\label{rmk:Z}
Note that $V$ is itself an $\a$-stable subordinator. (This will be relevant in our discussion on aging in Section~\ref{sec:age}
below---see Remark~\ref{rmk:age3} in that section.)
We may regard $V$ and $Z$ as processes in the random environment $\ups$. In fact, given $\ups$, they are both Markovian.
We may think of $\ups$ as the scaling limit of the (relevant) environment of the trap model. 
The overall distribution of $Z$ (integrated over the joint distribution of $\ups$ and $\{T_r;\,r\geq0\}$) makes it a self similar 
process of index $1$.
\end{rmk}

\begin{rmk}
 \label{rmk:const}
Let $c>0$ be the constant such that $\esp(\exp\{-\l \ups_1\})=\exp\{-c\l^\a\}$. One readily checks that the distribution
of $Z$ does not depend on that constant. 
Here and below we will denote also by $\P$ and $\esp$ the probability and expectation underlying the distributions of
$\ups$ and the $T_r$'s.
\end{rmk}

We are now ready to state our convergence 
result, but first some notation. For $\eps>0,\,t\geq0$, let 
\begin{equation}
\label{eq:ye}
\E^{(\eps)}_t=\eps q_\eps\E_{\eps^{-1}t}, 
\end{equation}
where $q_\eps$ is a slowly varying function at $0$ to be further specified below, 
and denote  $\E^{(\eps)}=(\E^{(\eps)}_t)$ and $Z=(Z_t)$. 
Let $D$, $D_T$ denote the spaces of c\`adl\`ag real functions defined on $[0,\infty)$, $[0,T]$ 
respectively. Let $d_T$ denote the $L_1$ distance in $D_T$, and $d=\sum_{n=1}^\infty2^{-n}(d_n\wedge1)$.

\begin{theo}
\label{teo:conv1}
Suppose Assumptions A and B are in force. 
Then as $\eps\searrow0$
\begin{equation}
\label{eq:conv1}
\E^{(\eps)}\to Z
\end{equation}
in distribution on $(D,d)$, where $(q_\eps)_{\eps>0}$ is nonincreasing, slowly varying at $0$ and satisfies
$1\leq q:=\lim_{\eps\to0}q_\eps\leq\infty$, with $q<\infty$ if and only if $X$ is transient.
\end{theo}

Notice that in the transient case we have linear scaling, and sublinear scaling in the recurrent case.
A more precise description of $q_\eps$ is given in~(\ref{eq:ae2}) below.

Here is a rough evaluation on the form of the spatial scaling in~(\ref{eq:ye}) leading to~(\ref{eq:conv1}). 
After $n$ steps, the random walk $X$ has explored $R_n\sim\rho_n$ new sites, each
one visited of the order of $1/r_n$ times, so the walk should be visiting a trap of depth of order $s_{\rho_n}=:v_n$, 
where $s_n$ is of the order of the largest of $n$ independent copies of $\tau_0$. Under our assumptions, the main 
contributions to the total time spent by the process to give $n$ steps, namely the times spent at the deepest 
traps encountered up to that number of steps, are roughly independent, so that total time should be of order $v_n/r_n=:\phi_n$. 
Then to times of order $\eps^{-1}$ correspond energies of order $v_{\nu_\eps}$, where
\begin{equation}\label{eq:nue}
\nu_\eps=\phi^{-1}(\eps^{-1}),
\end{equation}
and $\phi^{-1}$ is the inverse of $\phi$.
We thus should take
\begin{equation}\label{eq:ae1}
\frac1{v_{\nu_\eps}}=\frac1{r_{\nu_\eps}\phi_{\nu_\eps}}\sim\frac{\eps}{r_{\nu_\eps}}
\end{equation} 
as spatial scaling, where the latter aproximate equality follows from the inverse relation between $\phi$ and $\nu$.
From the assumptions on $r$ and the distribution of $\tau_0$, it follows that $\nu$ is nonincreasing, unbounded and regularly
varying at infinity. Since $r\in[0,1)$ is nonincreasing and slowly varying at infinity, we conclude that 
\begin{equation}\label{eq:ae2}
q_\eps=1/r_{\nu_\eps} 
\end{equation} 
has the stated properties.


The rough discussion of the last paragraph also gives rise to the form of the limiting process, once one realizes 
that (under our assumptions) the total times spent on the few contributing traps, when divided by the respective
means, and suitably scaled to account for the number of distinct visits, are approximately iid mean 1 exponentials.
A striking aspect is the universality of the limiting process even for cases where the scalings are distinct
(linear and sublinear as pointed out above).

A remark about the topology: small traps do not contribute to the limit, and were they to be completely disregarded 
the convergence would take place in the $J_1$/usual Skorohod topology, 
but they are there and mix up with large traps in a way which the $J_1$ topology 
(and other more usual ones, like the $M_1$ topology) is too fine to handle, and so we resort to a rougher topology.

One important technical aspect to consider in order to prove the scaling limit of an aging function like~(\ref{eq:af1}) 
is to show that for an arbitrary fixed time $t$, $\E^{(\eps)}_t$ converges to $Z_t$ as $\eps\to0$ in a strong enough form. 
This does not follow from Theorem~\ref{teo:conv1}, mainly due to the topological issue just discussed.
We state and prove such a convergence result, Lemma~\ref{lm:comp}, in Section~\ref{sec:age}, 
involving suitable versions of the relevant processes, 
before establishing the main results of that section, namely the aging results stated in 
Theorem~\ref{teo:age} (see also Remark~\ref{rmk:age4}). 
Our approach requires the strengthening of Assumptions A and B to transience of $X$. 
Scaling limit results for integrated versions of the aging functions herein considered, and potentially 
others, follow from Theorem 5, or rather the arguments in its proof, directly, under the original assumptions. 
See Subsection~\ref{ssec:iage}.

We close with a brief discussion on the stronger than annealed version for Theorem~\ref{teo:conv1}. Roughly, the annealed aspect of
Theorem~\ref{teo:conv1} follows from our approach of considering a version of the 
$\tau$ variables placed over the range of the underlying discrete random walk, thus fixing it, in such a way that they in a certain sense 
converge almost surely (to the increments of a stable subordinator, $\Upsilon$) --- see details in the proof of Theorem~\ref{teo:conv1} on
Section~\ref{sec:conv}. This turns out to be convenient for the analysis leading to Theorem~\ref{teo:conv1}, but when going back to the original
$\tau$'s, we only have an annealed result.
But of course, when we consider the distribution of the process given $\tau$, we integrate with respect to the underlying discrete random walk, 
and the averaging involved in this could lead to a stronger result. 
A condition is required, however --- annealed convergence is all we have in, say, the asymmetric simple one dimensional case. One is
introduced in Section~\ref{sec:str}, saying roughly that independent realizations of the trajectory of $X$ intersect little. With
that additional condition, we state
and prove stronger convergence results.

%% file: pre30.tex
\section{Preliminaries on random walks}
\label{sec:rw}
\setcounter{equation}{0}

In this section, we establish a few facts concerning discrete time random walks that follow from Assumptions  $A$ and  $B$ made in
Section~\ref{sec:mod}
above, as well as from other assumptions we will consider below. These results will be used later in the sections ahead.

Let $X=(X_n,\,n\geq0)$ be the random walk introduced in Section~\ref{sec:mod} above, and define 
\begin{equation}
\label{eq:run} 
u_n=\P(X_n=0),\quad U_n=\sum_{i=0}^nu_i,\quad L_n=\sum_{i=0}^n\1\{X_i=0\},\quad n\geq0.
\end{equation}
$L_n$ is the occupation time of the origin up to step $n$. We will also write $L_x$ for a positive real $x$, and it means $L_{\lf x\rf}$,
similarly for $U_x$ and $r_x$.

Our first remark is that 
\begin{equation}
\label{eq:ror} 
\rho_n=\sum_{k=0}^n r_k,
\end{equation}
with $\rho$ and $r$ defined in~(\ref{eq:range1},~\ref{rn}) above.
The formula is proved as follows (see~\cite{kn:S} page 36). First we have that
\begin{equation}\nn
R_n=\sum_{k=0}^n \1\{R_k=R_{k-1}+1\}\,, 
\end{equation}
where $R$ was defined in~(\ref{eq:range}) above. Upon noticing that 
\begin{equation}\nn
R_k=R_{k-1}+1\Leftrightarrow X_k-X_{k-1}\not=0 ; X_k-X_{k-2}\not=0,\ldots,X_k\not=0\,,
\end{equation}
we conclude that 
\begin{eqnarray*}
&\P(R_k=R_{k-1}+1)=&\\
&\P(\xi_k\not=0;\xi_k+\xi_{k-1}\not=0,\ldots,X_k\not=0)
=\P(X_1\not=0;X_2\not=0,\ldots,X_k\not=0)=r_k\,.&
\end{eqnarray*}
 It then follows from~(\ref{eq:ror}) and Assumption  $B$ that 
\begin{equation}
\label{eq:ronr} 
\lim_{n\to\infty}\frac{\rho_n}{nr_n}=1\,.
\end{equation}
(This readily follows from the first displayed equation on page 55 of~\cite{kn:Se}.)

Our second remark is the following result.

\begin{lm}
\label{lm:approx}
Under Assumption  $B$, and provided $\lim_{n\to\infty}r_n=0$ (that is, if $X$ is recurrent), the law of $r_n L_n$ approximates a mean 1 exponential
distribution as $n\to\infty$. 
\end{lm}

\noindent{\bf Proof} 

Given $u\geq0$, we have that
$$ \{r_n L_n>u\}=\{\eta_1+...+\eta_k\leq n\} $$ 
where the $\eta_j$, $j\geq1$, are the successive  increments of return times to the origin by $X$, and $k=\lf u/r_n\rf+1$. 

Therefore 
\begin{equation} \label{1}
\P(r_n L_n>u)=\P(\eta_1+...+\eta_k\leq n)=\P(\bar\eta_n\leq1),
\end{equation} 
where $\bar\eta_n=(\eta_1+...+\eta_k)/n$.

A straightforward computation of the Laplace transform of $\bar\eta_n$ yields
\begin{equation} \label{2}
\esp(e^{-\l\bar\eta_n})=\{\esp(e^{-\frac\l n\eta_1})\}^k=\left\{1-\frac\l n\int_0^\infty r(x)\,e^{-\frac\l n x}\,dx\right\}^k
=\left\{1-\int_0^\infty r(y n/\l)\,e^{-y}\,dy\right\}^k,
\end{equation} 
where  $r(x)=\P(\eta_1>x)=r_{\lf x\rf}$.

Assumption  $B$ and Theorems 2.6 and 2.7 of~\cite{kn:Se} imply that the integral on the right hand side of~(\ref{2}) is asymptotic
to $r_n$ as $n\to\infty$. This and the form of $k$ imply that
\begin{equation} \label{3}
\lim_{n\to\infty}\esp(e^{-\l\bar\eta_n})=e^{-u}
\end{equation} 
for all $\l>0$, and this implies that the law of $\bar\eta_n$ converges to that of a(n extended) random variable which takes the value $0$
with probability $e^{-u}$, and the value $\infty$ with the complementary probability. It follows that 
\begin{equation} \label{4}
\lim_{n\to\infty}\P(r_n L_n>u)=e^{-u}
\end{equation} 
for every $u>0$. $\square$

\begin{cor}
\label{cor:approx}
Under Assumption $B$, we have that
\begin{equation} \label{5}
\lim_{n\to\infty}r_n U_n=1.
\end{equation} 
\end{cor}

\noindent{\bf Proof}  

In the transient case, this follows from $r_n\to r_\infty>0$, $L_n\to L_\infty$, a Geometric random variable with mean $r_\infty^{-1}$, 
and monotone convergence.

In the recurrent case, from the first equality in~(\ref{1}), we find that
\begin{equation} \label{6}
\P(r_n L_n>u)\leq\P\!\left(\max_{1\leq j\leq k}\eta_j\leq n\!\right)=(1-r_n)^k.
\end{equation} 
From the form of $k$, and the fact that $\lim_{n\to\infty}(1-r_n)^{1/r_n}=e^{-1}$, we find that for every $c>e^{-1}$ and large enough $n$, the right
hand side of~(\ref{5}) is dominated by $c^u$ for all large enough $n$. Dominated convergence now yields
\begin{equation} \label{7}
r_n U_n=\int_0^\infty\P(r_n L_n>u)\,du\to\int_0^\infty e^{-u}\,du=1\,.
\end{equation} 
as $n\to\infty$, since, from~(\ref{eq:run}), $r_n U_n=r_n\esp(L_n)= \esp(r_nL_n)$. $\square$

\bigskip

Another corollary of Lemma~\ref{lm:approx} is as follows. For $x\in\Z^d$, let
\begin{equation} \label{71}
\ell(x,n)=\sum_{j=0}^n \1\{X_j=x\}\,T_j,
\end{equation}
where $T_1,T_2,\ldots$ are iid mean 1 exponential random variables independent of $X$.

\begin{cor}
\label{cor:approx1}
Under Assumption $B$, we have that $r_n\ell(0,n)$ converges weakly to a mean 1 exponential distribution.
\end{cor}

\noindent{\bf Proof}  

In the recurrent case, the result follows immediately from Lemma~\ref{lm:approx} and the law of large numbers, 
once we observe that 
\begin{equation} \label{8}
\ell(0,n)=\sum_{i=1}^{L_n} T'_j ,
\end{equation} 
where $T'_1,T'_2,\ldots$ are iid mean 1 exponential random variables independent of $X$.  

In the transient case, $L_n$ converges as $n\to\infty$ to a geometrically distributed random variable, say $L_\infty$. 
Also $\lim_{n\to\infty}r_n=r_\infty>0$, and one readily checks that $r_\infty\sum_{i=1}^{L_\infty} T'_j$ is a mean 1 exponential 
random variable. 
$\square$

\bigskip

\begin{lm}
\label{lm:cp}
Under Assumption $B$, we have that for every $0<a<b<\infty$
\begin{equation}
\label{eq:labn}
L_{bn}-L_{an}\to0
\end{equation}
in probability as $n\to\infty$.
\end{lm}

\begin{rmk}
 \label{rmk:labn}
Since $L_{bn}-L_{an}$ is an integer, we have that the probability of no return to $0$ of $X$ during $[an,bn]$ goes to 1 as $n\to\infty$ for every
fixed $0<a<b<\infty$.
\end{rmk}

\noindent{\bf Proof}

Let $N=N_n(a,b)$ denote the random variable on the 
 left of~(\ref{eq:labn}). Using the Markov property, one readily checks that the conditional
distribution of $L_{(b+1)n}-L_{an}$ given $N\geq1$ dominates the unconditional one of $L_n$. 
We conclude that
\begin{eqnarray}\nn
U_{(b+1)n}-U_{an}\=\esp[L_{(b+1)n}-L_{an}]\geq\esp[L_{(b+1)n}-L_{an};\,N\geq1]\\
\label{eq:cp1}
&\geq&\esp[L_n]\,\P(N\geq1)=U_n\,\P(N\geq1).
\end{eqnarray}

We then have that $\P(N\geq1)\leq(U_{(b+1)n}-U_{an})/U_n$ and the result follows from Assumption $B$ and~(\ref{5}).
$\square$

%% file: conv54.tex
\section{Convergence}
\label{sec:conv}

\setcounter{equation}{0}

This section is devoted to the proof of Theorem~\ref{teo:conv1}.
The main idea of the proof is to work with a 
version of the
process in which the environment is dynamically constructed along the
trajectory of the walk and coupled to a stable process. Once
this 
version is defined it is rather elementary to get the
convergence in a quenched form.

\begin{rmk}
 \label{rmk:inv}
We often in this and other sections consider the inverse $g$ of 
a (possibly random) monotonic unbounded function 
$f:{\mathbb D}\to[0,\infty)$, where ${\mathbb D}$ may be either 
$\N$ or $[0,\infty)$, defined (as usually)
\begin{equation}
  \label{eq:n}
  g(x)=\inf\{y\in{\mathbb D}:\,f(y)>x\}
\end{equation}
for $x\in[0,\infty)$. (Some times, we will write $f_n$ instead of $f(n)$.)
\end{rmk}

Let $X$ and $\tau$ be as in the previous sections. Assumptions $A$ and $B$ are in force.
In the build up to the proof of Theorem~\ref{teo:conv1}, we start by considering a particular construction of the law of $\X$ and $\E$
under $\P$, as follows.
Let $T_0,T_1,T_2,\ldots$ be a family of independent mean $1$ exponential random variables. Consider the following random function
$C:\N\to[0,\infty)$:
\begin{equation}
\label{eq:clock}
C_n=\sum_{i=0}^n\tau_{X_i} T_i,\quad n\geq0,
\end{equation}
and let $I$ denote its inverse. We may call $C=(C_n)$ the {\em clock process} associated to (this particular construction of) the trap model.
Now define for $t\geq0$
\begin{equation}
\label{eq:y}
Y_t=\tau_{X_{I_t}}.
\end{equation}
Note that $(X_{I_t})$  has the same law as $\X$ under $\P$, and $Y=(Y_t)$ has the same law as $\E$ under $\P$. Thus, making 
\begin{equation} 
\label{eq:yet}
Y^{(\eps)}_t=
\eps q_\eps Y_{\eps^{-1}t},\,t\geq0,
\end{equation}
we have that $Y^{(\eps)}=(Y^{(\eps)}_t)$ has the same law as $\E^{(\eps)}$ under $\P$.

We will work with a particular version of $Y^{(\eps)}$ where a specific version of the (scaled) $\tau$ random variables coupled to $\Upsilon$ are
effectively placed over the range of $X$.
We define this specific version in~(\ref{eq:bye}) below. The proof of Theorem~\ref{teo:conv1} is then divided into first showing that that is
indeed a version, and next establishing convergence of the version.
The ingredients of the version are the above defined $X$ and $\Upsilon$, and a family of iid mean 1 exponential random variables
\begin{equation}
 \{\tj_i,\,j\geq0,\,i\geq1\}.
\end{equation}

Let us start by enumerating the {\em full range of $X$}
\begin{equation}\label{eq:enum}
 \RR:=\cup_{n\geq1}\RR_n=:\{\txt_0,\txt_1,\ldots\},
\end{equation} 
in {\em chronological} order (in the natural way, i.e., given $x,y\in\RR$, we have $x<y$ iff $X$ hits $x$ before it does $y$);
let $\psi:\N\to\N$ be the map
\begin{equation}
\label{eq:phi}
\psi(n)=m\quad\mbox{iff}\quad X_n=\txt_m.
\end{equation}

We now consider properly scaled $\tau_0$-distributed random variables, to be eventually placed over $\RR$, following the order therein.
For that let us (re)introduce
\begin{equation}
  \label{eq:vnu}
v_n=s_{\rho_n}\quad\mbox{and}\quad \phi_n=\frac{v_n}{r_n},
\end{equation}
where
\begin{equation}
  \label{eq:sn}
  s_n=\inf\{t\geq0:Q(t,\infty)\leq n^{-1}\}.
\end{equation}

Notice that all of this functions (of $n$), namely $v,\,\phi,\,s$ are nondecreasing and unbounded. 
They were discussed above (in the paragraph of~(\ref{eq:ae1})); we have now a more precise definition of $s_n$. 
Due to the regularly varying characters of $Q$ (see~(\ref{eq:tail})) and $r$ (Assumption B), $v,\,\phi,\,s$ 
are all regularly varying with a common index $\a^{-1}$).

It follows from elementary properties of monotonicity and regular variation of the above functions that
\begin{eqnarray}
  \label{eq:nuinf}
\nu_\eps\to\infty
\end{eqnarray}
as $\eps\to0$, where $\nu_\eps$ was defined in~(\ref{eq:nue}) above.

Now let 
\begin{equation}
  \label{eq:te}
\te=1/\rho_{\nu_\eps},
\end{equation}
 and make, for $x\in\N$,
\begin{equation}
  \label{eq:taue}
  \tilde\tau_x^{(\eps)}=G^{-1}\!\left(\te^{-1/\a}(\ups_{\te x+\te}-\ups_{\te x})\right),
\end{equation}
where $G$ is defined by
\begin{equation}
  \label{eq:G}
  \P(\ups_1>G(y))=\P(\tau_0>y),\,y\geq0.
\end{equation}
By elementary properties of the subordinator $\ups$ and the definition of $G$, it readily follows that
$\{\tilde\tau_x^{(\eps)},\,x\in\N\}$ is an iid family with $\tilde\tau_0^{(\eps)}\stackrel{d}=\tau_0$, 
where ``$\stackrel{d}=$'' denotes equality in distribution.

We are now ready to define the version of the age process which will be used in the proof of Theorem~\ref{teo:conv1}. 
We start by defining the version of the clock process. For $n\geq0$, let
\begin{equation}
\label{eq:bclock}
\tilde C^{(\eps)}_n=\sum_{y\in\N} \tilde\tau_{y}^{(\eps)}
\sum_{i=1}^{L\!\left(\txt_{y},n\right)}\!\!\!\hat T^{(y)}_i,
\end{equation} 
and let $\tilde I^{(\eps)}$ denote the inverse of $\tilde C^{(\eps)}$. Then make
\begin{equation}
\label{eq:bye}
\tilde Y^{(\eps)}_t=
\tilde\tau^{(\eps)}_{\psi\left(\tilde I^{(\eps)}_t\right)},\,\,t\geq0.
\end{equation}

\begin{rmk}
 \label{rmk:expl}
One may think of the latter version as one in which the environmental variables of the age process
$\tilde\tau^{(\eps)}=\{\tilde\tau_{x}^{(\eps)},\,x\in\,\N\}$
are placed over the range of $X$, $\{\txt_{x},\,x\in\N\}$, respectively.
The coupled construction of $\tilde\tau^{(\eps)}$ yields its strong convergence, when properly rescaled,
to its limiting counterpart $\ups$.
\end{rmk}

\noindent{\bf Proof of Theorem~\ref{teo:conv1}}

The result follows readily from Propositions~\ref{prop:version} and~\ref{lm:conv} stated and proven below.

$\square$

\begin{prop}
\label{prop:version}
For every $\eps>0$, $\tilde Y^{(\eps)}$ and $Y$ have the same distribution.
\end{prop}

Let $\hat Y^{(\eps)}_t=\eps q_\eps\tilde Y^{(\eps)}_{\eps^{-1}t},$ $t\geq0$.

\begin{prop}
\label{lm:conv}
For almost every $\ups$
\begin{equation}
\label{eq:conve}
\hat Y^{(\eps)}\to Z
\end{equation}
as $\eps\searrow0$ in distribution on $(D,d)$.
\end{prop}

\begin{rmk}
 \label{rmk:aeups}
The distribution referred to in the statement of Proposition~\ref{lm:conv} is the joint one of $X$ and $\{\tj_i\}$ (with $\ups$ fixed).
\end{rmk}

\noindent{\bf Proof of Proposition~\ref{prop:version}}

As noted right below~(\ref{eq:G}) above, we have that $\{\tau_{\hat X_y};\,y\in\N\}$ is independent of $X$ and
\begin{equation}
\label{eq:seq}
\{\tau_{\hat X_y};\,y\in\N\}\stackrel{d}=\{\tilde\tau^{(\eps)}_{y};\,y\in\N\}.
\end{equation}
It thus follows that
\begin{equation}
\label{eq:clock1}
C_n\stackrel{d}=\sum_{z\in\Z^d}\tau_{z}\sum_{i=1}^{L(z,n)} T^{(z)}_i
=\sum_{y\in\N}\tau_{\hat X_y}\sum_{i=1}^{L(\hat X_y,n)} T^{(\hat X_y)}_i
\stackrel{d}=\sum_{y\in\N}\tilde\tau^{(\eps)}_{y}\sum_{i=1}^{L(\hat X_y,n)} T^{(\hat X_y)}_i
\stackrel{d}=\sum_{y\in\N}\tilde\tau^{(\eps)}_{y}\sum_{i=1}^{L(\hat X_y,n)}\hat T^{(y)}_i,
\end{equation}
as vectors indexed by $n$, where $\{T^{(z)}_i;\,z\in\Z^d, i\geq1\}$ is an iid family of mean one 
exponential random variables.
$\square$

\bigskip

\noindent{\bf Proof of Proposition~\ref{lm:conv}}

It will be implicit (and sometimes explicit) in the claims made below that they hold for a.e.-$\ups$.
We will use the symbols ``$P$`` and ``$E$`` to denote the probability measure and expectation associated to the
distribution of $X$ and $\{\tj_i\}$ -- referred to sometimes below as dynamical random variables.

We start by defining the set of {\em deep traps} or $\delta$-{\em traps}. 
For $x\geq0$, let $\mu(x)=\ups(x)-\ups(x-)$ and fix $\delta>0$ arbitrarily. Consider
\begin{equation}
\label{eq:taudelta}
\t_{\delta}=\{x\geq0:\, \mu(x)>\delta\}=\{x_1<x_2<\ldots\}.
\end{equation}
Let now $\yep_i=\lf\te^{-1}x_i\rf$, $i\geq1$, and define
\begin{equation}
\label{eq:tauedelta}
\teps_{\delta}=\{\yep_1,\yep_2,\ldots\}.
\end{equation}

\begin{figure}[htb]
\begin{center}
\input{yte.pspdftex}
\end{center}
\caption{Schematic realization of $\hat Y^{(\eps)}$. Largest values appear in full horizontal lines; other values are negligible as $\eps\to0$
and are not depicted.}
\label{yte}
\end{figure}
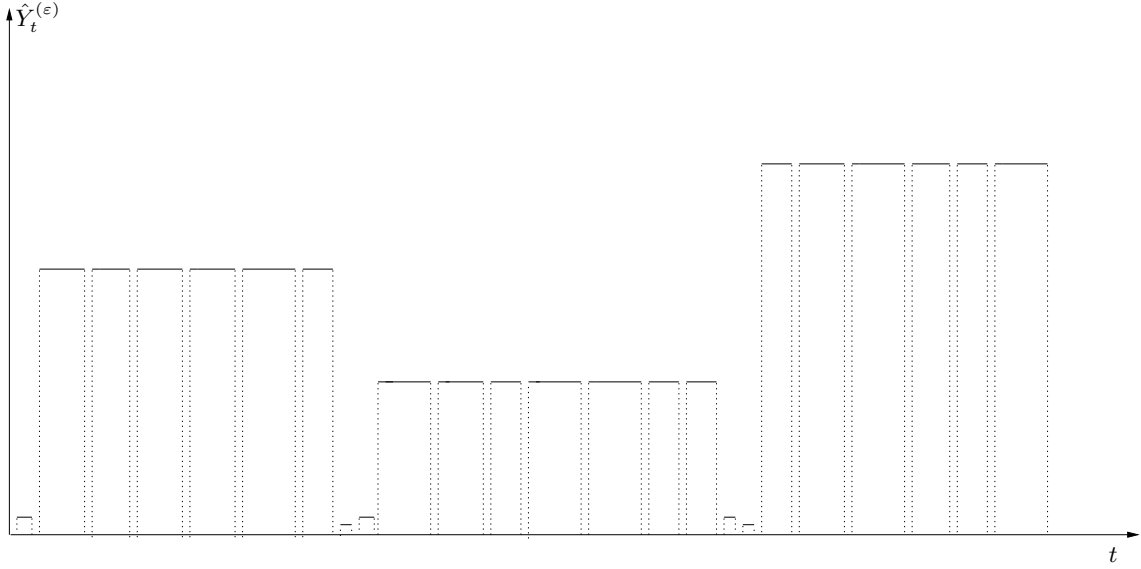

Our strategy will be to consider modified versions of $\hat Y^{(\eps)}$ and $Z$,
respectively $\hat Y^{(\eps,\d)}$ and $Z^{(\d)}$, where only deep traps contribute
(see~(\ref{eq:c5}) below), and then show on the one hand that a version of~(\ref{eq:conve})
holds for the modified processes (see Lemma~\ref{lm:conv1}), and on the other hand that 
the errors due to the replacing $\hat Y^{(\eps)}$ with $\hat Y^{(\eps,\d)}$ 
are negligible (see Lemma~\ref{lm:conv2} below).
One other aspect to be considered in adopting this modified version, besides the fact that the 
main contributions for $\hat Y^{(\eps)}$ come from its largest values, or $\delta$-traps, 
is that these contributions are interspersed with the lesser contributions of shallower traps. 
This is illustrated in Figure~\ref{yte}, where the $\d$-traps 
are represented in solid horizontal lines, and the heights of the lesser traps do not appear (but, 
perhaps inconsistently, the respective sojourn times do appear). This interspersion of small 
and large traps does not occur to $\hat Y^{(\eps,\d)}$ -- see Figures~\ref{tyed} and~\ref{tyedp}
below.

Let us first define a modified version of $\tilde Y^{(\eps)}$ in terms of a restricted clock process as follows. Let
\begin{equation}
\label{eq:rclock}
\tilde C^{(\eps,\d)}_n=\sum_{j\geq1} \tilde\tau_{\yep_j}^{(\eps)}\,\tilde T(\yep_j,n),
\end{equation}
where for $y\in\N$
\begin{equation}
\label{eq:tilt}
\tilde T(y,n)=\sum_{i=1}^{L\!\left(\txt_{y},n\right)}\!\!\!\hat T^{(y)}_i,
\end{equation}
and let $\tilde I^{(\eps,\d)}$ denote the inverse of $\tilde C^{(\eps,\d)}_n$ as a 
function of $n$. Then make 
\begin{equation}
\label{eq:byed}
\tilde Y^{(\eps,\d)}_t=
\tilde\tau^{(\eps)}_{\psi\left(\tilde I^{(\eps,\d)}_t\right)},\,\,t\geq0,
\end{equation}
and
\begin{equation}
\label{eq:hyed}
\hat Y^{(\eps,\d)}_t=
\eps q_\eps\tilde Y^{(\eps,\d)}_{\eps^{-1}t},\,\,t\geq0.
\end{equation}

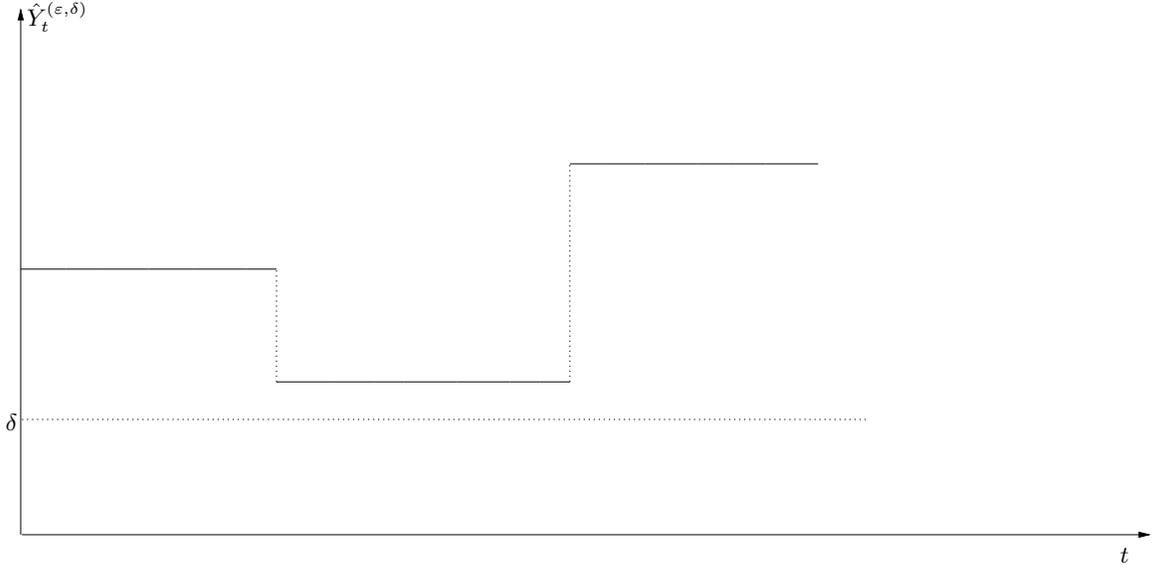
\begin{figure}[htb]
\begin{center}
\input{tyed.pspdftex}
\end{center}
\caption{Schematic realization of $\hat Y^{(\eps,\delta)}$. Only values of $\delta$-traps appear, with sojourn times corresponding to
actual times spent in those values by $\hat Y^{(\eps)}$.}
\label{tyed}
\end{figure}

\begin{figure}[htb]
\begin{center}
\input{tyedp.pspdftex}
\end{center}
\caption{Schematic realization of $\hat Y^{(\eps,\delta')}$, $\delta'<\delta$. See Figure~\ref{tyed}.}
\label{tyedp}
\end{figure}
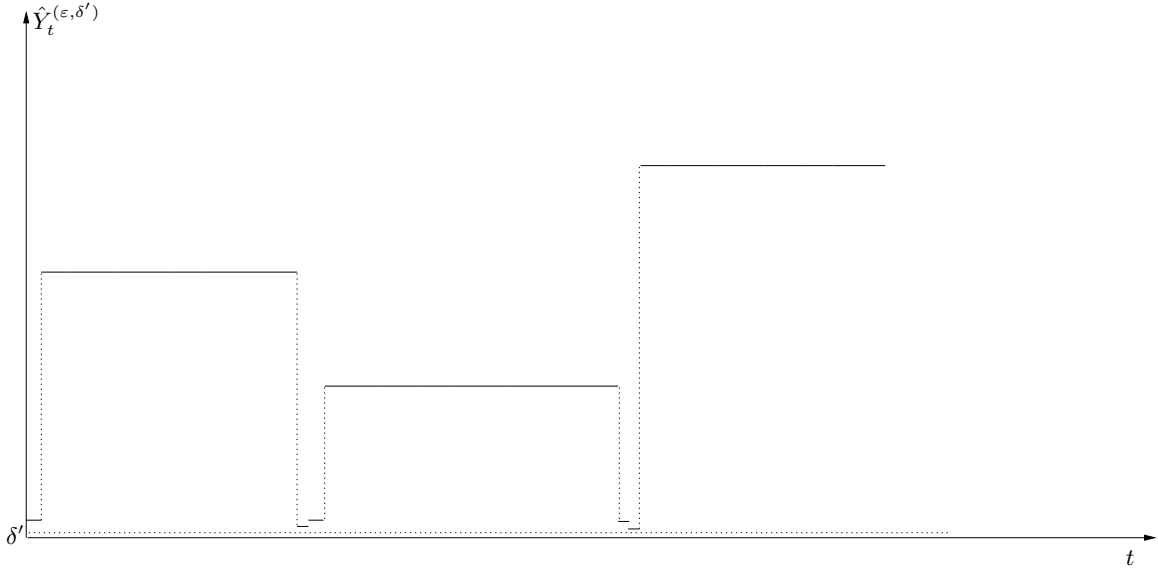


Let us now define the modified versions of $Z$. Let
\begin{equation}
\label{eq:c4}
\bvd_s=\sum_{i=1,2,\ldots:\atop{x_i\leq s}}\mu(x_i)\,T_{x_i},\quad s\geq0,
\end{equation}
and consider the inverse $\bvd$ as a function of $s$, $\bwd$. Then let
\begin{equation}
\label{eq:c5}
\bzd_t=\mu(\bwd_t),\quad t\geq0.
\end{equation}

Theorem~\ref{teo:conv1} readily follows from Lemmas~\ref{lm:conv1},~\ref{rmk:conv_bvd} 
and~\ref{lm:conv2} below.

$\square$


\begin{lm}
\label{lm:conv1}
Let $\delta>0$ be fixed.
Then as $\eps\searrow0$
\begin{equation}
\label{eq:conve1}
\byed\to\bzd
\end{equation}
in distribution on $(D,J_1)$, where $J_1$ is the usual Skorohod metric on $D$.
\end{lm}

\noindent{\bf Proof}

We start by fixing a positive integer $K$, and analysing the behavior of $\byed$ till it hits the $K+1$-st trap of
$\t_{\delta}$; in other words, till $X$ hits $\hat X_{\yep_{K+1}}$. 

For $j=1,2,\ldots$ let us define $\xe_j=\te\yep_j$ and 
\begin{equation}
\label{eq:mue}
\mue(\xe_j)=\eps q_\eps\tilde\tau^{(\eps)}_{\yep_j}.
\end{equation}
\noindent {\bf Claim 1} 
We claim that, outside an event of vanishing probability as $\eps\to0$, 
before $X$ hits $\hat X_{\yep_{K+1}}$, $\byed$ does nothing except 
\begin{itemize}
 \item[$i$.] visiting 
$\mue(\xe_1),\mue(\xe_2),\ldots\mue(\xe_K),$ 
successively, without backtracking;

 \item[$ii$.]  the sojourn times of those visits converge in distribution to an independent vector of $K$ 
exponential random variables with means $\mu(x_1),\ldots,\mu(x_K)$, respectively.
\end{itemize}

Proposition 3.1 of~\cite{kn:FIN} implies that for a.e.-$\ups$
\begin{equation}
\label{eq:mu}
\mue(\xe_j)\to\mu(x_j)\quad\mbox{as}\quad\eps\to0,
\end{equation}
for every $j\geq1$.
\begin{rmk}
 \label{rmk:fin}
Some matching of the notation of~\cite{kn:FIN} and the present one needs to be 
done in order to verify the above claim by resorting to that reference.
$V$ there corresponds to $\ups$ here. Respectively, $\eps$ corresponds to $\te$.
We have that $1/v_{\nu_\eps}(\sim\eps q_\eps)$ above corresponds to $c_{\te}$, 
with $c_\cdot$ introduced in (3.9) of~\cite{kn:FIN}, 
and our $\tilde\tau_x^{(\eps)}$ corresponds to $\tau_{x}^{(\te)}$ of that 
reference. (\ref{eq:mu}) follows then from the point process convergence
statement in Proposition 3.1 of~\cite{kn:FIN}.
\end{rmk}

Claim 1 above and~(\ref{eq:mu}) are the main ingredients of the proof of Lemma~\ref{lm:conv1}.

In order to justify Claim 1, let us introduce $\ze_i$, the hitting time of $\hat X_{\yep_i}$ 
by $X$, $i=1,2,\ldots$, and consider
\begin{equation}
\label{eq:l1}
L\!\left(\txt_{\yep_i},\ze_{i+1}\right),\,i=1,2,\ldots
\end{equation} 

Let us state a result concerning these hitting times, proven at the end of this section.
\begin{lm}
\label{lm:sclaim}
For all $s'>s>s''\geq0$ and $\d>0$, we have that outside an event of vanishing probability as $\eps\to0$
\begin{equation}
\label{eq:sc1}
L(\txt_{y},\,\nu_\eps s)=0 \mbox{ for all } y>\te^{-1}s',\,\mbox{ and }\,
L(\txt_{y},\,\nu_\eps s)>0 \mbox{ for all } y\in\N\cap[0,\te^{-1}s''].
\end{equation}
\end{lm}

\begin{rmk}
 \label{rmk:xitox}
The above statement is equivalent to
\begin{equation}
\label{eq:scalt}
\frac1{\nu_\eps} \ze_i\to x_i
\end{equation}
in probability as $\eps\to0$ for $i\geq1$.
\end{rmk}

It follows from Lemma~\ref{lm:sclaim} that, 
outside an event of vanishing probability as $\eps\to0$, 
$\bced_{\nu_\eps s}$ vanishes if $s<x_1$, and is restricted to the $k$ 
first terms if $x_k<s<x_{k+1}$, $k=1,\ldots,K$.
Thus, given arbitrary $r_1<s_1<\ldots<r_{K+1}<s_{K+1}$ such that $0<r_i<x_i<s_i$, $i=1,\ldots,K+1$, 
outside an event of vanishing probability as $\eps\to0$ 
we have that 
\begin{equation}
\label{eq:l2}
\nu_\eps r_i\leq\ze_{i}\leq \nu_\eps s_i,
\end{equation}
and thus for $i=1,\ldots,K$
\begin{equation}\label{eq:l3}
L\!\left(\txt_{\yep_i},\ze_{i}+\nu_\eps(r_{i+1}-s_i)\right)
\leq L\!\left(\txt_{\yep_i},\ze_{i+1}\right)
\leq L\!\left(\txt_{\yep_i},\ze_{i}+\nu_\eps(s_{i+1}-r_i)\right)\!.
\end{equation}

Furthermore, by Corollary~\ref{rmk:labn} above, we get the following.

{\noindent \bf Key facts} Outside an event of vanishing probability as $\eps\to0$, we have equalities
in~(\ref{eq:l3}); namely
\begin{equation}\label{eq:kf}
L\!\left(\txt_{\yep_i},\ze_{i}+\nu_\eps(r_{i+1}-s_i)\right)
= L\!\left(\txt_{\yep_i},\ze_{i+1}\right)
= L\!\left(\txt_{\yep_i},\ze_{i}+\nu_\eps(s_{i+1}-r_i)\right)\!,
\end{equation}
$i=1,\ldots,K$, and indeed also
\begin{equation}\label{eq:kf1}
L\!\left(\txt_{\yep_i},\ze_{i+1}\right)
=L\!\left(\txt_{\yep_i},\ze_{K+1}\right)\!,\, i=1,\ldots,K.
\end{equation}

The first part of the Claim 1 at the beginning of this proof is thus established.

To argue the second part, we start by observing that
the sojourn time of $\tilde Y^{(\eps,\d)}$ on $\yep_j$, $j=1,\ldots,K$, up until $X$ hits 
$\hat X_{\yep_{K+1}}$ is given by 
\begin{equation}
\label{bte}
\tilde\tau^{(\eps)}_{\yep_j}\sum_{i=1}^{L\!\left(\txt_{\yep_j},\,\ze_{K+1}\right)}\!\!\hat T^{(\yep_j)}_{i},
\end{equation}
so the respective time spent by $\hat Y^{(\eps,\d)}$ is $\eps$ times that, which can be rewritten as
\begin{equation}
\label{bte1}
\mue(\xe_j)\,\tilde T^{(\eps)}_j,
\end{equation}
where
\begin{equation}
\label{bte2}
\tilde T^{(\eps)}_j=\,r_{\nu_\eps}\!\!\!\!\!\sum_{i=1}^{L\!\left(\txt_{\yep_j},\,\ze_{K+1}\right)}\!\!\hat T^{(\yep_j)}_{i}.
\end{equation}
(Let us recall that $q_\eps^{-1}=r_{\nu_\eps}$.)

From the above we conclude that outside an event of vanishing probability as $\eps\to0$, 
up until $X$ hits $\hat X_{\yep_{K+1}}$, we have that
\begin{equation}
\label{eq:tzed}
\byed_t=\bzed_t:=\mue(\bwed_t),
\end{equation}
where $\bwed$ is the inverse of 
\begin{equation}
\label{eq:l7}
\bved_s=\sum_{j=1,2,\ldots:\atop{x_j\leq s}}\mue(\xe_j)\,\bte_j.
\end{equation}

We now argue that $\bte_j,\,j=1,2,\ldots,$ are asymptotically 
independent mean 1 exponential random variables as $\eps\to0$.
By the key fact~(\ref{eq:kf}) and Remark~\ref{rmk:labn} above, and the Markov property,
we have that for $a>0$, outside an event of vanishing probability as $\eps\to0$, 
each $\bte_j$ coincides with
\begin{equation}
r_{\nu_\eps}\!\!\!\!\sum_{i=1}^{L\!\left(\txt_{\yep_j},\,\ze_{j}+\nu_\eps(r_{j+1}-s_j)\right)}\!\!\hat T^{(\yep_j)}_{i}
=r_{\nu_\eps}\!\!\!\!\sum_{i=1}^{L\!\left(\txt_{\yep_j},\,\ze_{j}+a\nu_\eps\right)}\!\!\hat T^{(\yep_j)}_{i},
\end{equation}
and each of the latter quantities is equally distributed with
\begin{equation}
\label{eq:l5}
r_{\nu_\eps}\!\!\sum_{i=1}^{L\!\left(0,a\nu_\eps\right)}\!\!\!\hat T^{(1)}_i,
\end{equation}
$j=1,\ldots,K$.

Corollary~\ref{cor:approx1} 
now implies that for any fixed $a>0$ the distribution of~(\ref{eq:l5})
converges to a mean $1$ exponential one as $\eps\to0$. This and~(\ref{eq:mu}) 
establish the exponentiality part of Claim 1.$ii$.

To establish the independence part, it is enough to argue that for $j\geq2$ we have that $\bte_j$
is asymptotically independent of the vector 
\begin{equation}
\label{eq:l7a}
(\bte_i,\,i=1,\ldots, j-1).
\end{equation}
By the key facts~(\ref{eq:kf}-\ref{eq:kf1}), we have that,
outside an event of vanishing probability as $\eps\to0$,
\begin{equation}
\bte_j=
r_{\nu_\eps}\!\!\!\!\sum_{i=1}^{L\!\left(\txt_{\yep_j},\,\ze_{j}+\nu_\eps(r_{j+1}-s_j)\right)}\!\!\hat T^{(\yep_j)}_{i}.
\end{equation}
Since $\ze_{j}$ is a stopping time, the strong Markov property yields the independence
of the latter random variable and the vector in~(\ref{eq:l7a}), and the claimed asymptotic
independence is established.
Claim 1 is thus established.

In order to complete the argument for this proof, we first observe that we may replace $D$ 
with $D_T$, $T>0$ arbitrary. 

Note that from the above arguments, it follows that for every $s\geq0$ fixed, we have that
\begin{equation}\label{vtov}
 \bved_s\to V^{(\d)}_s
\end{equation}
in distribution as $\eps\to0$.
It is clear that 
\begin{equation}
\label{eq:vdinf}
\bvd_s\to\infty \mbox{ as } s\to\infty 
\end{equation}
almost surely, so given $\eta>0$, we may choose $S\notin\t_\d$ such that for all small enough $\eps$
\begin{equation}
\label{eq:st}
P(\bved_S>T)>1-\eta. 
\end{equation}
Let $K=|\t_{\delta}\cap[0,S]|$. Then outside an event of probability at
most $\eta$ the distance $D_T$ is bounded above by $D_{\bved_S}$, and combining this with Claim 1 at the beginning
of this proof, 
we get that outside an event of probability at most $2\eta$, uniformly in $\eps$ small, within $[0,T]$,
$\byed$ visits $\mue(\xe_1),\ldots,\mue(\xe_K)$ successively without backtracking (not necessarily all of them
on $[0,T]$),
with sojourn times converging to independent exponential random variables
with means $\mu(x_1),\ldots,\mu(x_K)$, respectively. 
The result follows.
$\square$

\begin{rmk}
 \label{rmk:conv2}\mbox{ }
\begin{enumerate}
\item We claimed above that the convergence in~(\ref{vtov}) holds for fixed times $s$.
In the next section we will need a stronger version of that convergence, namely one that holds for 
the trajectories, under the $J_1$ metric. It is pretty clear from the ingredients at end of the above proof that
\begin{equation}
\label{conv4}
\lim_{\eps\to0}(\bved_t)\stackrel{d}=(\bvd_t)
\end{equation}
on $(D,J_1)$. Let us sketch a brief argument. Indeed, $(\bved_t)$ is a jump process with jumps located at $\xe_1,\xe_2,\ldots$ with respective sizes
$\mue(\xe_1)\bte_1,\mue(\xe_2)\bte_2,\ldots$. As established above, the jump
locations converge to $x_1,x_2,\ldots$, and the jump sizes converge in distribution to independent exponentials of respective
means $\mu(x_1),\mu(x_2),\ldots$, which is a description of $(\bvd_t)$, and the convergences of jump locations and sizes
imply $J_1$-convergence.

\item Another point to be used below: it immediately follows from~(\ref{eq:l2}) that if $s\notin\t_\d$ then
\begin{equation}
\label{cv}
\eps\bced_{\nu_\eps s}=\bved_s
\end{equation} 
outside an event of vanishing probability as $\eps\to0$. It follows from~(\ref{vtov}) that,
for every fixed $s\geq0$, $\eps\bced_{\nu_\eps s}\to V^{(\d)}_s$
in distribution as $\eps\to0$.

In particular, we may conclude that given $T,\eta,\delta>0$, there exist $\eps_0,S>0$
such that $\P(\eps\bce_{\nu_\eps S}\leq T)\leq\P(\eps\bced_{\nu_\eps S}\leq T)\leq\eta$ for all $\eps<\eps_0$.

\end{enumerate}

\end{rmk}

\begin{lm}
 \label{rmk:conv_bvd}
We have that
\begin{equation}
\label{eq:c5a}
(\bvd_t)\to(V_t)\, \mbox{ and }\, (\bzd_t)\to(Z_t)
\end{equation}
almost surely on $(D,J_1)$ as $\delta\to0$. 
\end{lm}

\noindent{\bf Proof}

Using the fact that 
$\sum_{x\in\t_{\delta}^c\cap[0,T]}\mu(x)\to0$ as $\delta\to0$ for arbitrary $T$,
the proof follows readily.

$\square$

We will claim several times below that for certain random variables depending on two parameters $\eps$ and $\d$,
called generically now $\Xi^{(\eps,\d)}$, it holds that $\lim_{\delta\to0}\limsup_{\eps\to0}\Xi^{(\eps,\d)}=0$
in probability. This means that for every $\eta>0$, we have that
$\lim_{\delta\to0}\limsup_{\eps\to0}\P(\Xi^{(\eps,\d)}>\eta)=0$.

\begin{lm}
\label{lm:conv2}
\begin{equation}
\label{eq:c6}
\lim_{\delta\to0}\limsup_{\eps\to0}d((\byed_t),(\hye_t))=0
\end{equation}
in probability. 
\end{lm}

\noindent{\bf Proof}

To establish~(\ref{eq:c6}), we introduce a further auxiliary process. The reason for this is the following
difficulty in the direct comparison between $\byed$ and $\hye$. The trajectories of those two processes are
somewhat poorly {\em aligned} when the latter one is visiting $\d$-traps --- which is the only occasion they 
can agree, since the former process lives on $\d$-traps. By poor alignment of these processes we mean that 
there are times when $\hye$ is visiting
a given $\d$-trap, while $\byed$ is 
visiting a different $\d$-trap.
This is perhaps apparent on a comparison of Figures~\ref{yte} and~\ref{tyed}.
With the aim of minimizing this {\em bad} misalignment, we introduce a process,
$\vyed$, living on $\d$-traps which is {\em better} aligned with $\hye$ in the sense that the poor alignment 
described above does not take place, and which is easily compared to $\byed$ as well. More on this discussion 
after the definition of $\vyed$ next.

Let  $\bar x_0=0$, $\bar x_i=(x_i+x_{i+1})/2$, $i\geq1$, and define
\begin{equation}
\label{eq:vyed}
\vyed_t=\mue(\xe_i),\,\mbox{ if }\, \eps\bce_{\nu_\eps\bar x_{i-1}}<t\leq\eps\bce_{\nu_\eps\bar x_{i}},\,\,i\geq1.
\end{equation}

\begin{figure}[htb]
\begin{center}
\input{yted.pspdftex}
\end{center}
\caption{Schematic realization of $\bar Y^{(\eps,\delta)}$. Only values of $\delta$-traps appear, with sojourn times (slightly) longer than
actual times spent in those values by $\hat Y^{(\eps)}$. See Figure~\ref{tyed}.}
\label{byed}
\end{figure}
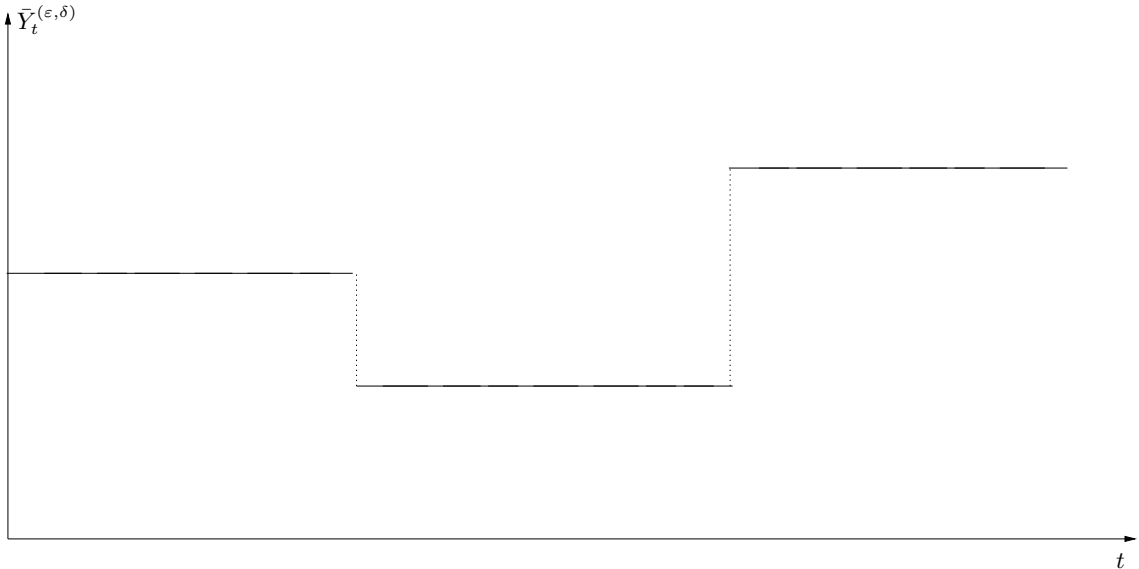

\begin{figure}[htb]
\begin{center}
\input{ytedp.pspdftex}
\end{center}
\caption{Schematic realization of $\bar Y^{(\eps,\delta')}$. See Figure~\ref{tyedp}.}
\label{byedp}
\end{figure}
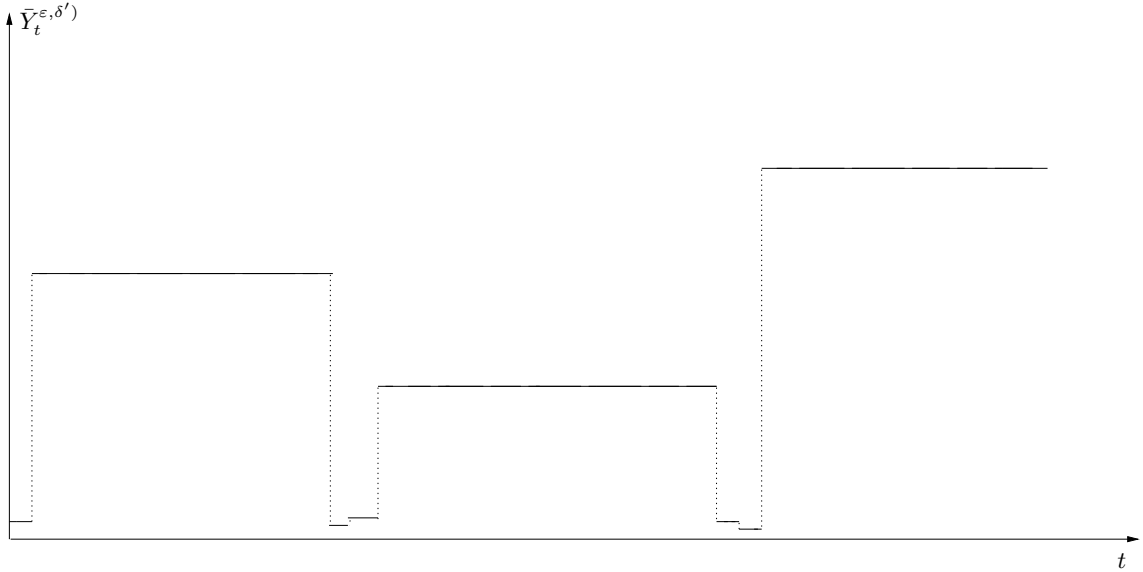

Coming back to the alignment issue, we note that on any given finite time interval, outside an event of vanishing 
probability as $\eps\to0$,
everytime $\hye$ is visiting a $\d$-trap, $\vyed$ is visiting the same trap. Since $\vyed$ lives on
$\d$-traps, the $d$ distance between $\vyed$ and $\hye$ on a given time interval is, outside an event of vanishing 
probability as $\eps\to0$, bounded above by the size of the deepest $\d$-trap in that interval multiplied by the total 
time spent by $\hye$ outside $\d$-traps during that interval.

The comparison between $\vyed$ and $\byed$ is also quite simple, perhaps simpler. Both processes live on $\d$-traps, which 
they visit in the same order, with distinct sojourn times of order 1, which approach each other in the limit as $\eps\to0$. 
So we indeed have the vanishing of the $J_1$ distance of $\vyed$ and $\byed$ in the limit as $\eps\to0$.

We make these arguments more precisely now. We have that~(\ref{eq:c6}) follows from
\begin{eqnarray}
\label{eq:c7}
&\lim_{\delta\to0}\limsup_{\eps\to0}d((\hye_t),(\vyed_t))=0,&\\
\label{eq:c8}
&\lim_{\delta\to0}\limsup_{\eps\to0}J_1((\vyed_t),(\byed_t))=0,&
\end{eqnarray}
in probability.

To establish~(\ref{eq:c7}),  we claim that it is enough to consider $d_T$ instead of $d$, with
\begin{equation}
\label{eq:N}
T=T(\eps,S)=\eps\sum_{y\in\N} \tilde\tau_{y}^{(\eps)}\,\tilde T(y,\nu_\eps S),
\end{equation}
$S$ fixed arbitrarily.  To justify the claim, it suffices to argue that $T\to\infty$ in probability as $S\to\infty$ 
uniformly in $\eps$ at a
neighborhood of the origin. This can be done as follows. $T$ clearly dominates $\eps\bced_{\nu_\eps S}$ 
(see~(\ref{eq:rclock}) above); it then follows
from~(\ref{cv}) and~(\ref{eq:vdinf}) that
the latter quantity diverges in probability as $S\to\infty$ uniformly in $\eps$ around the origin. 
This closes the argument for the claim. 

Reasoning now again as in the proof of Lemma~\ref{lm:conv1}, we find that outside 
an event of vanishing probability as $\eps\to0$, $\hye$ and $\vyed$ coincide within $[0,T]$ 
whenever the first process is in a $\d$-trap. Since the set of $\d$-traps visited during $[0,T]$ 
is contained in $[0,S+1]$, we conclude that the $d_T$ distance between the two processes is bounded 
above by 
\begin{equation}
\label{eq:c9}
\max\{\mue(\xe_i),\,i\geq1,\,\xe_i\leq S+1\}\,\,\eps\!\!\!\!
\sum_{y\in\N\setminus\teps_\delta} \tilde\tau_{y}^{(\eps)}\,\bte(y,\nu_\eps S).
\end{equation}
(\ref{eq:c7}) then follows from~Lemma~\ref{lm:neg} below and the fact that 
the max is over a bounded set uniformly in $\eps$ and~(\ref{eq:mu}).

To establish~(\ref{eq:c8}), we again replace $D$ by $D_T$, this time $T$ deterministic, but otherwise 
arbitrarily fixed. For an arbitrary $\eta>0$, let $S$ be as in the second point of Remark~\ref{rmk:conv2}. 
Then, arguing as in the proof of Lemma~\ref{lm:conv1} (see the Claim 1 at the beginning of that 
proof, and also the paragraph of~(\ref{eq:tzed})), on the event that $\eps\bced_{\nu_\eps S}>T$ and outside 
an event of vanishing probability as $\eps\to0$, $(\byed_t)_{t\in[0,T]}$ successively visits the set of states 
$\{\mue(\xe_i),\,i\geq1:\,\xe_i\leq S+1\}$ (not necessarily all of them by time $T$, but in any case in that order),
with respective sojourn times 
$\{\mue(\xe_i)\,\bte_i,\,i\geq1:\,\xe_i\leq S+1\}$.
The same is of course true of $(\vyed_t)$, except that the sojourn times are given by 
$\{\tilde S^{(\eps)}_i:=\eps\bce_{\nu_\eps\bar x_{i}}-\eps\bce_{\nu_\eps\bar x_{i-1}},\,i\geq1:\,\xe_i\leq S+1\}$.                    
Furthermore, for $i\geq1$ such that $\xe_i\leq S+1$, we have that outside an event whose probability vanishes as $\eps\to0$,
$\tilde S^{(\eps)}_i\geq\hat S^{(\eps)}_i:=\mue(\xe_i)\, \bte_i$ for such $i$, 
and for the same $i$ the difference between $\tilde S^{(\eps)}_i$ and $\hat S^{(\eps)}_i$ is bounded above by

\begin{equation}
\label{eq:c10}
\eps\!\!\!\!\sum_{x\in\te\N\cap[\bar x_{i-1},\,\bar x_{i}]\setminus\{\xe_i\}}
\tilde\tau_{\te^{-1}x}^{(\eps)} \,\bte(\te^{-1}x,\nu_\eps S)
\leq
\eps\!\sum_{y\in\N\setminus\teps_\delta}\tilde\tau_{y}^{(\eps)}\,\bte(y,\nu_\eps S),
\end{equation}
and the latter expression vanishes in probability as $\eps\to0$
by Lemma~\ref{lm:neg}. The result follows since $\eta$ is arbitrary.
$\square$

\begin{rmk}
 \label{rmk:byed}
It follows from arguments in the last paragraph of the above proof together with other arguments that for all $t>0$ fixed
\begin{equation}
\label{eq:comp1a}
\lim_{\delta\to0}\limsup_{\eps\to0}P(\byed_t\ne\vyed_t)=0.
\end{equation}
Indeed for $t<T$, with $T$, $\eta$ and $S$ as in the proof of~(\ref{eq:c8}) above (see last paragraph of the above proof),
and setting $K=\max\{i\geq1:\,x_i<S+2\}$, let 
\begin{equation}
\tse_j=\sum_{i=1}^j\tilde S^{(\eps)}_i,\quad\hse_j=\sum_{i=1}^j\hat S^{(\eps)}_i,\,j=1,\ldots,K,
\end{equation}
where $\tilde S^{(\eps)}$ and $\hat S^{(\eps)}$ were defined right above~(\ref{eq:c10}). Then, since every random variable involved
is continuous and $\tilde S^{(\eps)}_j\geq\hat S^{(\eps)}_j$, $1\leq j\leq K$ (outside an event whose probability vanishes as $\eps\to0$),
we have that $\{\byed_t\ne\vyed_t\}$ 
is almost surely contained in 
$\cup_{j=1}^K\{\hse_{j}<t<\tse_{j}\}$. Now the fact argued at end of the above proof that the difference 
between $\hat S^{(\eps)}_{j}$ and $\tilde S^{(\eps)}_{j}$ vanishes as $\eps\to0$ in probability for $j=1,\ldots,K$,
and the fact that $\hat S^{(\eps)}_1,\ldots,\hat S^{(\eps)}_K$ are asymptotically independent, lead readily to the completion of the argument.
\end{rmk}

\begin{lm}
\label{lm:neg}
We have that
\begin{equation}
\label{eq:c2a}
\lim_{\delta\to0}\limsup_{\eps\to0}\,\,\eps\!\!\!\!\sum_{y\in\N\setminus\teps_\delta} 
\tilde\tau_{y}^{(\eps)}\,\bte(y,\nu_\eps s)=0
\end{equation}
in probability for every $s\geq0$. 
\end{lm}

\noindent{\bf Proof}

We start by resorting to the Markov property and Corollary~\ref{cor:approx} 
to find that for every $s>0$
\begin{equation}
\label{eq:expcurly}
r_{\nu_\eps}E(\bte(y,\nu_\eps s))\leq r_{\nu_\eps}E(L(0,\nu_\eps s))=
r_{\nu_\eps}U_{\nu_\eps s}\to1\quad\mbox{as}\quad\eps\to0.
\end{equation}

We next resort to Lemma~\ref{lm:sclaim} to restrict the sum in~(\ref{eq:c2a}) on $x\leq \te^{-1}s'$. 
By~(\ref{eq:expcurly}) the expectation of the restricted sum is bounded above by constant times 
\begin{equation}
\label{eq:c3}
\sum_{y\in\N\cap(\teps_\delta)^c\cap[0,\te^{-1}s']}\!\!\!\!\eps q_\eps\tilde\tau_{y}^{(\eps)} 
\end{equation}
(recall that $q_\eps=r_{\nu_\eps}^{-1}$). We now claim that the $\lim_{\delta\to0}\limsup_{\eps\to0}$ 
of~(\ref{eq:c3}) vanishes almost surely. In order to use results of~\cite{kn:FIN} for that, let us
extend to $\te\N$ the measure $\mue$ defined above for rescaled $\d$-trap sites only (see~(\ref{eq:mue})).
For $x\in\te\N$, let
\begin{equation}
 \label{eq:muee}
\mue(x)=\eps q_\eps\tilde\tau_{\te^{-1}x}^{(\eps)}.
\end{equation} 
Here we abuse notation by writing $\mue(\cdot)$, $\mu(\cdot)$ instead of $\mue(\{\cdot\})$, $\mu(\{\cdot\})$.

We may then rewrite the sum~(\ref{eq:c3}) as 
\begin{equation}
\label{eq:sf}
\mue([0,s'])-\mue([0,s']\setminus\te\t_\d^{(\eps)})
\end{equation}
Now by~(\ref{eq:mu}) we have that the second term of~(\ref{eq:sf}) converges to
$\mu([0,s']\cap\t_\d)$ almost surely as $\eps\to0$. It follows from Proposition 3.1 of~\cite{kn:FIN} (the vague 
convergence part) that $\mue([0,s'])\to\mu([0,s'])=\ups(s')$ almost surely as $\eps\to0$. See Remark~\ref{rmk:fin} 
above. Here we may use the fact that deterministic points like $s'$ are continuity points of $\ups$ almost surely.

We then have that~(\ref{eq:c3}) converges almost surely to 
\begin{equation}
\label{eq:sf1}
\mu([0,s']\setminus\t_\d) 
\end{equation}
as $\eps\to0$. It is again a standard fact that the expression in~(\ref{eq:sf1}) vanishes 
almost surely as $\d\to0$, and the claim and lemma follow. 
$\square$

\medskip

\noindent{\bf Proof of Lemma~\ref{lm:sclaim}}

We start by pointing out that
\begin{equation}
\label{eq:condl}
L(\txt_y,\nu_\eps s)=0 \,\mbox{ if and only if } y>R_{\nu_\eps s}. 
\end{equation}
Given $\eta>0$ to be chosen below, consider the event
\begin{equation}
\label{eq:event}
\left\{\left|\frac{R_{\nu_\eps s}}{\rho_{\nu_\eps s}}-1\right|>\eta\right\}.
\end{equation}
Assumption A implies that the probability of this event vanishes as $\eps\to0$.

Now write $\rho_{\nu_\eps s}$
as
\begin{equation}
\label{eq:rest1}
\frac{\rho_{\nu_\eps s}}{\nu_\eps s\,r_{\nu_\eps s}}\,
s\,
\frac{r_{\nu_\eps s}}{r_{\nu_\eps}}\,
\frac{\nu_\eps r_{\nu_\eps}}{\rho_{\nu_\eps}}\,\te^{-1}.
\end{equation}
From~(\ref{eq:ronr}) and~(\ref{eq:nuinf}), we get that the $\lim_{\eps\to0}$ of the first and third quotients in~(\ref{eq:rest1}) both equal 1.
Assumption B 
implies that the $\lim_{\eps\to0}$ of the second quotient equals 1.

Now let us choose $\eta>0$ small enough so that $(1+2\eta)s\leq s'$ and $(1-2\eta)s\geq s''$. 
From the conclusion of the previous paragraph, we get that, for all $\eps$ small enough, the events
\begin{eqnarray}
\label{eq:events}
&\{R_{\nu_\eps s}\geq\te^{-1} x\,\mbox{ for some }\,x>s'\}=\{R_{\nu_\eps s}>\te^{-1}s'\}&\\
&\{R_{\nu_\eps s}\leq\te^{-1} x\,\mbox{ for some }\,x\leq s''\}=\{R_{\nu_\eps s}\leq\te^{-1}s''\}&
\end{eqnarray}
are contained in event~(\ref{eq:event}), which, as already pointed out, has vanishing probability as $\eps\to0$. 

Since outside the event 
\begin{equation}
\label{eq:union}
\{R_{\nu_\eps s}\geq\te^{-1} x\,\mbox{ for some }\,x>s'\}\cup\{R_{\nu_\eps s}\leq\te^{-1} x\,\mbox{ for some }\,x\leq s''\}
\end{equation}
we have that
(\ref{eq:condl}) yields~(\ref{eq:sc1}), the result follows. 
$\square$

%% file: yte.pspdftex
\begin{picture}(0,0)%
\includegraphics{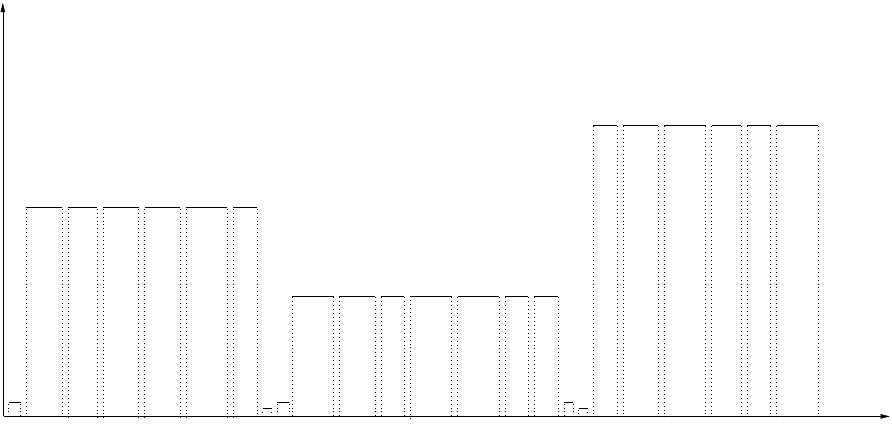}%
\end{picture}%
\setlength{\unitlength}{2072sp}%
\begingroup\makeatletter\ifx\SetFigFont\undefined%
\gdef\SetFigFont#1#2#3#4#5{%
  \reset@font\fontsize{#1}{#2pt}%
  \fontfamily{#3}\fontseries{#4}\fontshape{#5}%
  \selectfont}%
\fi\endgroup%
\begin{picture}(13569,6738)(859,-6358)
\put(14056,-6271){\makebox(0,0)[lb]{\smash{{\SetFigFont{9}{10.8}{\familydefault}{\mddefault}{\updefault}{\color[rgb]{0,0,0}$t$}%
}}}}
\put(976,149){\makebox(0,0)[lb]{\smash{{\SetFigFont{9}{10.8}{\familydefault}{\mddefault}{\updefault}{\color[rgb]{0,0,0}$\hat Y^{(\varepsilon)}_t$ }%
}}}}
\end{picture}%

%% file: tyed.pspdftex
\begin{picture}(0,0)%
\includegraphics{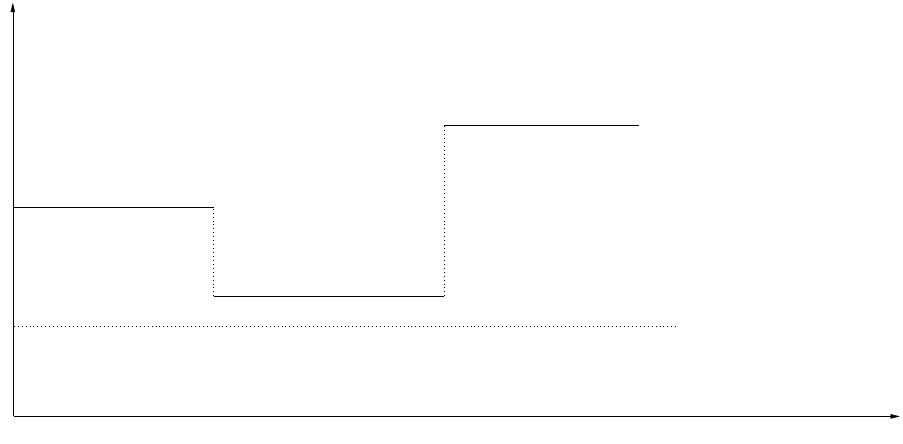}%
\end{picture}%
\setlength{\unitlength}{2072sp}%
\begingroup\makeatletter\ifx\SetFigFont\undefined%
\gdef\SetFigFont#1#2#3#4#5{%
  \reset@font\fontsize{#1}{#2pt}%
  \fontfamily{#3}\fontseries{#4}\fontshape{#5}%
  \selectfont}%
\fi\endgroup%
\begin{picture}(13722,6738)(706,-6358)
\put(14056,-6271){\makebox(0,0)[lb]{\smash{{\SetFigFont{9}{10.8}{\familydefault}{\mddefault}{\updefault}{\color[rgb]{0,0,0}$t$}%
}}}}
\put(721,-4696){\makebox(0,0)[lb]{\smash{{\SetFigFont{9}{10.8}{\familydefault}{\mddefault}{\updefault}{\color[rgb]{0,0,0}$\delta$}%
}}}}
\put(976,149){\makebox(0,0)[lb]{\smash{{\SetFigFont{9}{10.8}{\familydefault}{\mddefault}{\updefault}{\color[rgb]{0,0,0}$\hat Y^{(\varepsilon,\delta)}_t$ }%
}}}}
\end{picture}%

%% file: tyedp.pspdftex
\begin{picture}(0,0)%
\includegraphics{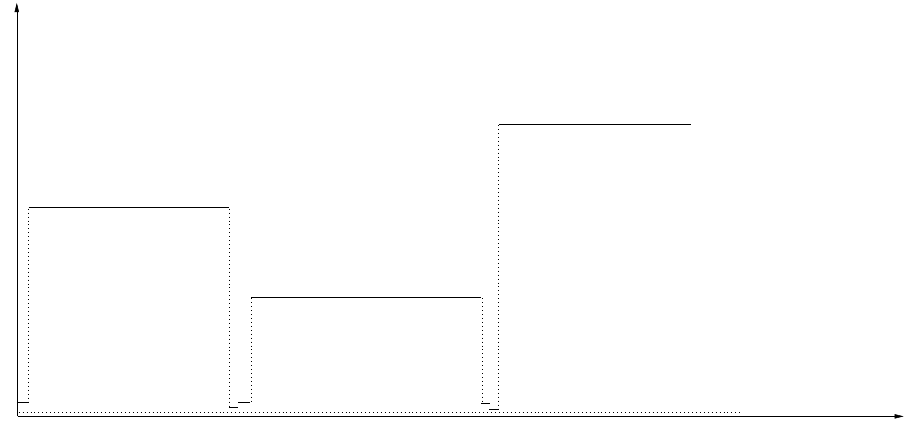}%
\end{picture}%
\setlength{\unitlength}{2072sp}%
\begingroup\makeatletter\ifx\SetFigFont\undefined%
\gdef\SetFigFont#1#2#3#4#5{%
  \reset@font\fontsize{#1}{#2pt}%
  \fontfamily{#3}\fontseries{#4}\fontshape{#5}%
  \selectfont}%
\fi\endgroup%
\begin{picture}(13782,6738)(646,-6358)
\put(14056,-6271){\makebox(0,0)[lb]{\smash{{\SetFigFont{9}{10.8}{\familydefault}{\mddefault}{\updefault}{\color[rgb]{0,0,0}$t$}%
}}}}
\put(661,-6016){\makebox(0,0)[lb]{\smash{{\SetFigFont{9}{10.8}{\familydefault}{\mddefault}{\updefault}{\color[rgb]{0,0,0}$\delta'$ }%
}}}}
\put(976,149){\makebox(0,0)[lb]{\smash{{\SetFigFont{9}{10.8}{\familydefault}{\mddefault}{\updefault}{\color[rgb]{0,0,0}$\hat Y^{(\varepsilon,\delta')}_t$ }%
}}}}
\end{picture}%

%% file: yted.pspdftex
\begin{picture}(0,0)%
\includegraphics{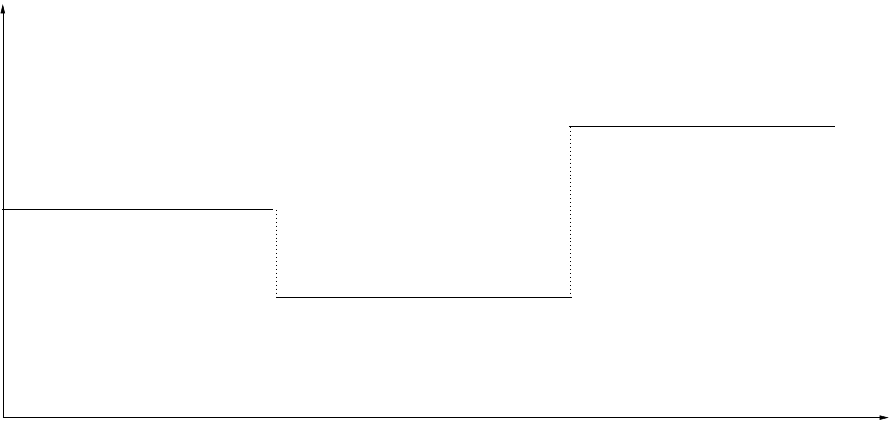}%
\end{picture}%
\setlength{\unitlength}{2072sp}%
\begingroup\makeatletter\ifx\SetFigFont\undefined%
\gdef\SetFigFont#1#2#3#4#5{%
  \reset@font\fontsize{#1}{#2pt}%
  \fontfamily{#3}\fontseries{#4}\fontshape{#5}%
  \selectfont}%
\fi\endgroup%
\begin{picture}(13554,6768)(874,-6373)
\put(14176,-6286){\makebox(0,0)[lb]{\smash{{\SetFigFont{9}{10.8}{\familydefault}{\mddefault}{\updefault}{\color[rgb]{0,0,0}$t$}%
}}}}
\put(1021,164){\makebox(0,0)[lb]{\smash{{\SetFigFont{9}{10.8}{\familydefault}{\mddefault}{\updefault}{\color[rgb]{0,0,0}$\bar Y^{(\varepsilon,\delta)}_t$}%
}}}}
\end{picture}%

%% file: ytedp.pspdftex
\begin{picture}(0,0)%
\includegraphics{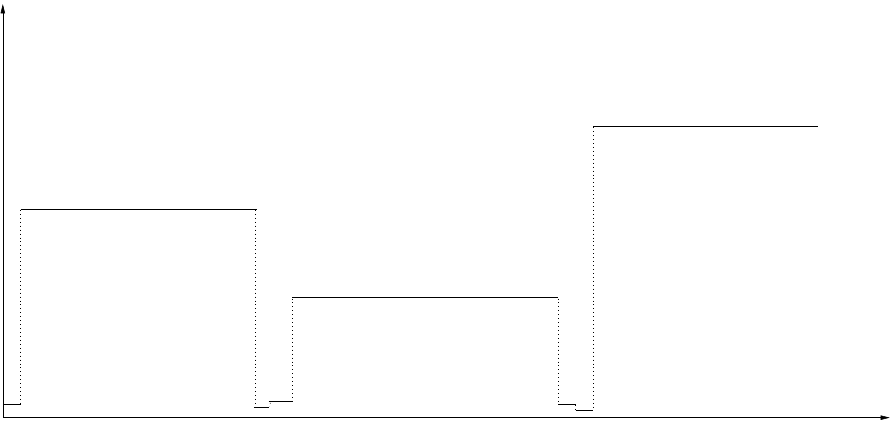}%
\end{picture}%
\setlength{\unitlength}{2072sp}%
\begingroup\makeatletter\ifx\SetFigFont\undefined%
\gdef\SetFigFont#1#2#3#4#5{%
  \reset@font\fontsize{#1}{#2pt}%
  \fontfamily{#3}\fontseries{#4}\fontshape{#5}%
  \selectfont}%
\fi\endgroup%
\begin{picture}(13569,6768)(859,-6373)
\put(14161,-6286){\makebox(0,0)[lb]{\smash{{\SetFigFont{9}{10.8}{\familydefault}{\mddefault}{\updefault}{\color[rgb]{0,0,0}$t$}%
}}}}
\put(1021,164){\makebox(0,0)[lb]{\smash{{\SetFigFont{9}{10.8}{\familydefault}{\mddefault}{\updefault}{\color[rgb]{0,0,0}$\bar Y^{\varepsilon,\delta')}_t$}%
}}}}
\end{picture}%

%% file: age54.tex
\section{Aging}
\label{sec:age}
\setcounter{equation}{0}

We will consider the following two aging functions of $(Y_t)$.
\begin{eqnarray}
\label{eq:br}
&\bar R(s,t)=\P(Y_t=Y_{t+s}),&\\
\label{eq:bpi}
&\bar \Pi(s,t)=\P(Y_t=Y_{t+r}\,\mbox{ for all }r\in[0,s]),&
\end{eqnarray}
with the following result.
\begin{theo}
\label{teo:age}
If $X$ is transient, then there exist non-trivial functions $R,\Pi:[0,\infty)\to(0,1]$ such that
\begin{eqnarray}
\label{eq:age_r}
&\lim_{t\to\infty}\bar R(\theta t,t)=R(\theta),&\\
\label{eq:age_pi}
&\lim_{t\to\infty}\bar \Pi(\theta t/q_{1/t},t)=\Pi(\theta).&
\end{eqnarray}
\end{theo}

\begin{rmk}
\label{rmk:trans}
We do not have any strong reason to believe that the result does not hold generally under Assumptions A and B.
The restriction to transient processes is technical. Our argument below requires, roughly speaking, that during 
time $t$ a {\em single} visit of  $Y$ to a deep trap lasts a length of time of order $t$, and this occurs only 
in the transient case. In this case, both assumptions A and B hold, with $\rho_n\sim n/q$, $r_n\sim 1/q$, and 
$q_{1/t}\sim q$, where $1<q<\infty$ is the same constant appearing in the statement of Theorem~\ref{teo:conv1}.

We state below a (weaker) version of this result for integrated versions of the aging functions, or equivalently,
for those functions looked at suitable random times.
\end{rmk}

\begin{rmk}
\label{rmk:age1}
Let $\tilde X_t=X_{I_t}$. In the literature one has rather considered the aging functions
\begin{eqnarray}
\label{eq:r}
&R(s,t)=\P(\tilde X_t=\tilde X_{t+s}),&\\
\label{eq:pi}
&\Pi(s,t)=\P(\tilde X_t=\tilde X_{t+r}\,\mbox{ for all }r\in[0,s]).&
\end{eqnarray}
In case $\tau_0$ is a continuous random variable, then we of course have the identities $R(\cdot,\cdot)=\bar R(\cdot,\cdot)$ 
and $\Pi(\cdot,\cdot)=\bar\Pi(\cdot,\cdot)$, but not otherwise. In any case, one can show that aging results like~(\ref{eq:age_r},\ref{eq:age_pi}) 
hold for $R(\cdot,\cdot)$ and $\Pi(\cdot,\cdot)$ as well, with $R(\cdot)$ and $\Pi(\cdot)$ as limiting aging functions, respectively.
\end{rmk}

\begin{rmk}
\label{rmk:age2}
$R$ and $\Pi$ turn out to be identical. We have
\begin{equation}
\label{eq:pir} 
R(\theta)=\Pi(\theta)=\frac{\sin(\pi\alpha)}{\pi}\int_{\theta/(1+\theta)}^1s^{-\a}(1-s)^{\a-1}\,ds.
\end{equation}
See~(\ref{eq:age2}, \ref{eq:h2}) and Remark~\ref{rmk:age3} below. 
\end{rmk}

In order to prove Theorem~\ref{teo:age}, we will naturally consider the rescaled version $\hye$ of $Y$ with the special strongly converging rescaled
environment (see~(\ref{eq:bye}) above). One ingredient of the proof of~(\ref{eq:age_r}) is a comparison to $\vyed$ (see~(\ref{eq:vyed}) above) as
follows.

\begin{lm}
\label{lm:comp}
For all $t>0$ fixed, if $X$ is transient, then we have that
\begin{equation}
\label{eq:comp1}
\lim_{\delta\to0}\limsup_{\eps\to0}P(\hye_t\ne\vyed_t)=0
\end{equation}
for a.e.~$\ups$.
\end{lm}

\noindent{\bf Proof}

As in previous arguments, we will leave implicit many times below that claims hold for a.e.~$\ups$.
$P$ and $E$ remain as notation for the probability and expectation with respect to the distribution
of the dynamical random variables ($X$ and the $\hat T_i^{(j)}$'s).
We first consider for $T>0$
\begin{equation}
\label{eq:comp2}
\int_0^TP(\hye_s\ne\vyed_s)\,ds=E\int_0^T\1\{\hye_s\ne\vyed_s\}\,ds.
\end{equation}

We will argue as in the proofs of Lemmas~\ref{lm:conv1} and~\ref{lm:conv2} above. Let us first fix an arbitrary $\eta>0$, and then choose 
$S$ as in the second point of Remark~\ref{rmk:conv2}. Then on the event that $\eps\bce_{\nu_\eps S}>T$ and outside an event of vanishing 
probability as $\eps\to0$, the integral on the right of~(\ref{eq:comp2}) is bounded above by
\begin{equation}
\label{eq:comp3}
\sum_{x\in\te\N\cap(\teps_\delta)^c\atop x\leq S+1}\mue(x) \,\bte(\te^{-1}x,\nu_\eps S).
\end{equation}
Using now~(\ref{eq:expcurly}) and the fact argued below~(\ref{eq:c3}), we have that the $\lim_{\delta\to0}\limsup_{\eps\to0}$
of the expectation of the integral is bounded above by $T\eta$. Since $\eta$ is arbitrary, we have that
\begin{equation}
\label{eq:comp4}
\lim_{\delta\to0}\limsup_{\eps\to0}\int_0^TP(\hye_s\ne\vyed_s)\,ds=0.
\end{equation}

We now fix $\delta'>\delta$ and define $I=\min\{i\geq1:\,x_i\in\t_{\delta'}\}$.  
Let $\ze$ be as in the proof of Lemma~\ref{lm:conv1} (see paragraph of~(\ref{eq:l1})). 
Let $t>0$ be fixed and condition on 
\begin{equation}
\label{eq:barc}
\bar C^{(\eps)}:=\eps\bce_{\zeta^{(\eps)}_{I+1}-1}=\mue(x_I^{(\eps)})\hat T(\yep_I,\ze_I)+\Delta^{(\eps)},
\end{equation} 
where the latter summand, $\Delta^{(\eps)}$, is defined by this equality. 
Note that the former summand is an exponential
random variable of mean $\mue(x_I^{(\eps)}) r_{\nu(\eps^{-1})}$, and that
given $X$  (and $\ups$) the summands are 
absolutely continuous random variables independent of each other. 

It thus follows that
\begin{eqnarray}\nn
P(\hye_t\ne\vyed_t|X)&\leq&
\int_0^t\left(\int_0^sb_\eps e^{-b_\eps(s-r)}f^{(\eps)}(r)\,dr\right)P(\hye_{t}\ne\vyed_{t}|\bar C^{(\eps)}=s,X)\,ds\\
\label{eq:comp5}
&+&P(\bar C^{(\eps)}\geq t|X),
\end{eqnarray}
where $b_\eps^{-1}=\mue(x_I^{(\eps)}) r_{\nu(\eps^{-1})}$ and $f^{(\eps)}$ is the density of $\Delta^{(\eps)}$ given $X$. 
The above integral is thus upper bounded by
\begin{equation}
\label{eq:comp6}
b_\eps\int_0^tP(\hye_{t}\ne\vyed_{t}|\bar C^{(\eps)}=s,X)\,ds.
\end{equation}
Now, by the independence of increments of $(\bce_n)$ given $X$, the probability inside the latter integral can be written as 
\begin{equation}
\label{eq:comp7}
P(\ddye_{t-s}\ne\bryed_{t-s}|X),
\end{equation}
where $(\ddye_t)$ and $(\bryed_t)$ are defined as $(\hye_{t})$ and $(\vyed_{t})$ respectively with 
$(\bce_{n+\xi^{(\eps)}_{I+1}-1}-\bce_{\xi^{(\eps)}_{I+1}-1})$ replacing $(\bce_n)$. We thus obtain 
(after integrating in $X$)
\begin{equation}
\label{eq:comp8} 
\lim_{\delta\to0}\limsup_{\eps\to0}\int_0^tP(\ddye_{t-s}\ne\bryed_{t-s})\,ds
=\lim_{\delta\to0}\limsup_{\eps\to0}\int_0^tP(\ddye_{s}\ne\bryed_{s})\,ds=0
\end{equation}
similarly as we did~(\ref{eq:comp4}). 
And we have that $\lim_{\delta\to0}\liminf_{\eps\to0}b_\eps^{-1}\geq \delta'/q>0$, 
since $x_I\in\t_{\delta'}$ and $\lim_{\eps\to0}r_{\nu(\eps^{-1})}=1/q>0$. 
We thus get that the 
$\lim_{\delta\to0}\limsup_{\eps\to0}$ of the 
first term in~(\ref{eq:comp5}) vanishes, and since $\delta'$ is arbitrary and 
$\lim_{\delta'\to0}\limsup_{\eps\to0}\bar C^{(\eps)}=0$ in probability, 
as can be readily checked, the result follows.
$\square$

\begin{rmk}
\label{rmk:trans2}
In the recurrent case $b_\eps$ is not bounded, and thus the above argument breaks down.
\end{rmk}


\noindent{\bf Proof of Theorem~\ref{teo:age}}

We will first replace $\hye$ by $\byed$ and then resort to Lemma~\ref{lm:comp} and Remark~\ref{rmk:byed}.
By Lemma~\ref{lm:conv1},
we have that
\begin{equation}
\label{eq:age1} 
\lim_{\delta\to0}\lim_{\eps\to0}\byed=Z
\end{equation}
in distribution on $(D,J_1)$. From that and Lemma~\ref{lm:comp} and Remark~\ref{rmk:byed} above, we claim that
\begin{equation}
\label{eq:age2} 
\lim_{\eps\to0}\bar R(\theta\eps^{-1},\eps^{-1})=\lim_{\delta\to0}\lim_{\eps\to0}\P(\byed_1=\byed_{1+\theta})=\P(Z_1=Z_{1+\theta}):=R(\theta).
\end{equation}
(\ref{eq:age_r}) then follows.  
The only point of the claim that needs arguing is the second equality. 
We first point out that from the construction of $\byed$ (see~(\ref{eq:c5})), since the $\mu(x_i)$'s are almost surely all distinct,
we have from~(\ref{eq:mu}) that, for all fixed $\delta$ and all small enough $\eps$, the probability in the second
term in~(\ref{eq:age2}) equals
\begin{equation}
\label{eq:age2a} 
\P(\byed_1=\byed_{1+r}\,\mbox{ for all }r\in[0,\theta])
\end{equation}
plus a small error, and the latter probability equals
\begin{equation}
\label{eq:age2aa} 
\P([1,1+\theta]\cap\mbox{range of }\bved=\emptyset).
\end{equation}

It readily follows from~(\ref{conv4}) and Remark~\ref{rmk:conv2}.1
that the second term in~(\ref{eq:age2}) equals
\begin{equation}
\label{eq:age2c} 
\P([1,1+\theta]\cap\mbox{range of }V=\emptyset)=\P(Z_1=Z_{1+r}\,\mbox{ for all }r\in[0,\theta]).
\end{equation}
Now the right hand side of~(\ref{eq:age2c}) equals that of~(\ref{eq:age2}) (since the $\mu(x_i)$'s are almost surely all distinct).
(\ref{eq:age_r}) is then settled.

In the above argument, we felt the need to go through~(\ref{eq:age2a}) and~(\ref{eq:age2aa}), since the (indicators of the) events in the 
the second probability in~(\ref{eq:age2}) and in~(\ref{eq:age2a}) are not almost surely continuous on $(D,J_1)$, but so is the event
in~(\ref{eq:age2aa}).

As regards~(\ref{eq:age_pi}), we have that
\begin{equation}
\label{eq:hpe}
\bar \Pi(\theta\eps^{-1}/q_\eps,\eps^{-1})=\P(Y_{\eps^{-1}}=Y_{\eps^{-1}+r}\,\mbox{ for all }r\in[0,\theta\eps^{-1}/q_\eps])
=\esp\left(e^{-\theta/(\eps q_\eps Y_{\eps^{-1}})}\right)=\esp\left(e^{-\theta/\hye_1}\right),
\end{equation}
and Lemma~\ref{lm:comp} implies that
\begin{equation}
\label{eq:h2}
\lim_{\eps\to0}\esp\left(e^{-\theta/\hye_1}\right)=\lim_{\d\to0}\lim_{\eps\to0}\esp\left(e^{-\theta/\vyed_1}\right)
=\esp\left(e^{-\theta/Z_1}\right)=:\Pi(\theta),
\end{equation}
where the second equality follows from~(\ref{eq:age1}), which implies marginal convergence
in distribution (since each fixed
deterministic time is almost surely a continuity point of $Z$).
$\square$

\begin{rmk}
\label{rmk:age3}
 One quickly checks that
\begin{equation}
\label{eq:age3} 
\P(Z_1=Z_{1+\theta})=\P(Z_1=Z_{1+r}\,\mbox{ for all }r\in[0,\theta]),
\end{equation}
which equals $\P([1,1+\theta]\cap\mbox{ range of }V=\emptyset)$, as noted in~~(\ref{eq:age2c}) above.
One can then obtain the right hand side of~(\ref{eq:pir})
as an expresssion for the latter probability. (This may be readily seen to follow from Proposition 3.1 in~\cite{kn:B}, 
since $V$ is an $\a$-stable subordinator; see Remark~\ref{rmk:Z} above.)
We further notice that we can write
\begin{equation}
\label{eq:age4} 
\P(Z_1=Z_{1+\theta})=\esp\left(e^{-\theta/Z_1}\right).
\end{equation}
See Remark~\ref{rmk:age2} above.  (\ref{eq:age3}) and~(\ref{eq:age4}) give us the Laplace transform of $1/Z_1$. An expression for the density of
that variable can be found in (5.97) of~\cite{kn:FM}.
\end{rmk}

\begin{rmk}
\label{rmk:age4}
Another aging function which is natural on one side and not as considered in the literature as the above ones on the other side, and also 
fits well in the above picture, is the following one.
\begin{equation}
\label{eq:q}
\textstyle \Omega(s,t)=\P(\sup_{r\in[0,t]}Y_r<\sup_{r\in[0,t+s]}Y_{r})
\end{equation}
It was suggested in~\cite{kn:FIN} as a ``measure of the prospects for novelty in the system``. 
Since $\sup_{r\in[0,t]}\hye_r=\sup_{r\in[0,t]}\vyed_r$ if $\hye_t=\vyed_t$, from Lemma~\ref{lm:comp},~(\ref{eq:c8}), Lemma~\ref{lm:conv1}
and~(\ref{eq:c5a}),
we have that
\begin{equation}
\label{eq:qage}
\textstyle \lim_{t\to\infty}\Omega(\theta t,t)=\P(\sup_{r\in[0,1]}Z_r<\sup_{r\in[0,1+\theta]}Z_{r})=:\Omega(\theta).
\end{equation}
This is an example where the limiting aging function requires full use of the process $Z$; in the previous cases, 
the limits could be expressed in terms of the (clock) process $V$ alone. 
We could not find an explicit expression for the right hand side of~(\ref{eq:qage}).
\end{rmk}


\subsection{Integrated aging results}
\label{ssec:iage}

By considering integrated aging functions, or equivalently aging functions looked at random times, we may circumvent
the difficulties exhibited by recurrence --- see Remarks~\ref{rmk:trans} and~\ref{rmk:trans2}
above. As an example, let us consider one such integrated aging function, and state the corresponding result.

Let
\begin{equation}
\label{eq:ir}
\bar{\mathfrak R}(\lambda,\mu)=\esp[\bar R(\lambda{\mathcal T},\mu {\mathcal T})]=\int_0^{\infty}e^{-t}\bar R(\lambda t,\mu t)\,dt,
\end{equation}
with $\bar R$ as in~(\ref{eq:br}) above, and ${\mathcal T}$ a mean one exponential random variable independent of every other
random variable in the problem.

\begin{theo}
\label{teo:iage}
If $X$ satisfies Assumptions  $A$ and $B$, then 
\begin{equation}
\label{eq:iage}
\lim_{\lambda,\mu\to\infty\atop{\lambda/\mu\to\theta}}\bar{\mathfrak R}(\lambda,\mu)=R(\theta),
\end{equation}
with $R$ as in Theorem~\ref{teo:age} above.
\end{theo}

\noindent{\bf Proof}

Let us for simplicity take $\lambda=\theta\mu$.
We may then write
\begin{equation}
\label{eq:ir1}
\bar{\mathfrak R}(\lambda,\mu)
=\esp\int_0^{\infty}e^{-t}\,{\mathbf 1}\{Y_{\mu t}=Y_{(\mu+\lambda) t}\}\,dt
=\esp\int_0^{\infty}e^{-t}\,{\mathbf 1}\{\hye_t=\hye_{(1+\theta)t}\}\,dt,
\end{equation}
where $\eps=\mu^{-1}$.

Given $0<\eta<1$,~(\ref{eq:c7}) allows us to bound the right hand side of~(\ref{eq:ir1}) above and below by $(1\pm\eta)$ times
\begin{equation}
\label{eq:ir2}
\lim_{\eps\to0}\esp\int_0^{\infty}e^{-t}\,{\mathbf 1}\{\byed_t=\byed_{(1+\theta)t}\}\,dt
=\esp\int_0^{\infty}e^{-t}\,{\mathbf 1}\{Z^{(\delta)}_t=Z^{(\delta)}_{(1+\theta)t}\}\,dt,
\end{equation}
respectively, as soon as $\delta$ is close enough to $0$, where the latter equality follows from Lemma~\ref{lm:conv1}. 
Since $\eta$ is arbitrary, by Lemma~\ref{rmk:conv_bvd} the left hand side of~(\ref{eq:ir1}) equals
\begin{equation}
\label{eq:ir3}
\esp\int_0^{\infty}e^{-t}\,{\mathbf 1}\{Z_t=Z_{(1+\theta)t}\}\,dt=\int_0^{\infty}e^{-t}\,\P(Z_t=Z_{(1+\theta)t})\,dt,
\end{equation}
since every fixed $t$ is almost surely a continuity point of $Z$. Now by the self similarity of index 1 exhibited by $Z$, 
the probability inside the integral does not depend on $t>0$, and the result follows, since
$\P(Z_1=Z_{(1+\theta)})=R(\theta)$. $\square$

\medskip

Similar results can be argued for integrated versions of the aging functions $\bar\Pi$ and $\Omega$ above.

%% file: str53.tex
\section{Stronger convergence}
\label{sec:str}

\setcounter{equation}{0}

In this section, we strengthen the convergence results of Section~\ref{sec:conv} under an additional condition, 
which we now explain.

Let $X'=(X'_n)_{n\geq0}$ a random walk independent from and equally distributed with $X$ and define 
\begin{equation}
\label{eq:int} 
\I_n=\I_n(X,X')=\{z\in\Z^d:\,X_i=X'_j=z\mbox{ for some }0\leq i,j\leq n\}=\RR_n(X)\cap\RR_n(X'),\, n\geq0,
\end{equation}
as the set of intersection points of the paths of $X$ and $X'$ up to (discrete) time $n$ (it can be seen also as indicated as the intersection of the
ranges of $X$ and $X'$ up to time $n$). Let now
\begin{equation}
\label{eq:nint} 
I_n=|\I_n|
\end{equation}
be the number of such intersection points. The additional condition we impose, in order that the results
of this section hold, is as follows.

\noindent{\bf Assumption} $\mathbf C$
\begin{equation}
\label{eq:adcond} 
\frac{I_n}{\esp(R_n)}\to0\mbox{ in probability as }n\to\infty.
\end{equation}

\begin{rmk}
 \label{rmk:adcond}
The expectation of the quotient in~(\ref{eq:adcond}) can be reexpressed as
\begin{equation}
\label{eq:adcond1} 
\frac{\sum_{x\in\Z^d}[\P(T_x\leq n)]^2}{\sum_{x\in\Z^d}\P(T_x\leq n)},
\end{equation}
where $T_x=\inf\{n\geq0:X_n=x\}$. We readily find that~(\ref{eq:adcond}) holds in either the general $d\geq2$ transient case, 
or the one dimensional  non integrable increment, transient case,
since in both these cases $\limsup_{\|x\|\to\infty}\P(T_x<\infty)=0$ (see e.g.~Proposition 25.3 in~\cite{kn:S}), and in general
$\lim_{n\to\infty}\esp(R_n)=\infty$.

It also holds for two dimensional mean zero, finite second moment random walks from results in~\cite{kn:L}.
We are uncertain about other recurrent planar walks, as well as 1-stable 1-dimensional recurrent walks. 
(See Remark~\ref{rmk:lln} above.)
\end{rmk}

Let $\B_u$ be the class of bounded uniformly continuous real functions on $(D,d)$. Here is the main result of this section.
Let $Y^{(\eps)}$ and $Z$ be as in Theorem~\ref{teo:conv1} above.

\begin{theo}
\label{teo:str}
Under Assumptions $A$, $B$ and $C$, for every $F\in\B_u$, we have
\begin{equation}
\label{eq:str}
\esp\left[\left.F(Y^{(\eps)})\right|\tau\right]\to\esp\left[F(Z)\right],
\end{equation}
in probability as $\eps\to0$.
\end{theo}

\begin{rmk}
\label{rmk:nec_cond}
As anticipated at the end of Section~\ref{sec:mod}, a condition like~(\ref{eq:adcond}) is needed for the validity of the above result. A case
where Assumptions $A$ and $B$ are satisfied, but not Assumption $C$, and~(\ref{eq:str}) does not hold, is when $X$ is one dimensional simple
asymmetric. This is particularly clear in the totally asymmetric case, when $C_n$ (see~(\ref{eq:clock})) is a partial sum of i.i.d.~random variables
in the basin of attraction of an $\a$-stable law, $\a\in(0,1)$, in which case it is well known to only converge, when properly rescaled, in
distribution. This prevents a result of the form of~(\ref{eq:str}) in that case.
\end{rmk}

\noindent{\bf Proof of Theorem~\ref{teo:str}}

Let $F\in\B_u$ be fixed. We may and will restrict to $F$ with bounded support, say $[0,T]$, where $T>0$ is arbitrary.

It follows from Theorem~\ref{teo:conv1} that 
\begin{equation}
\label{eq:s1}
\esp\left[F(Y^{(\eps)})\right]\to\esp\left[F(Z)\right].
\end{equation}
We will use this and~(\ref{eq:adcond}) to get~(\ref{eq:str}).

The strategy is to construct two sets of versions of $Y^{(\eps)}$. In each set of versions,
the different versions have independent dynamical random variables. The distinction is on the
environmental variables. On one set of versions the respective environments are also independent
among distinct versions, 
so that the versions are fully independent of one another: for this set the empirical 
mean of $F$ over the different versions yields the left hand side of~(\ref{eq:s1}) in the 
limit as the number of versions grows. 
On the other set of versions, we have a single environment for all the versions, so the empirical 
mean of $F$ over the different versions yields the left hand side of~(\ref{eq:str}).
However, this single environment 
can be constructed in a coupled way to the independent environments of the 
first set of versions, 
so that the difference between
the empirical mean over the first set of versions and that over the second set of versions
vanishes as $\eps\to0$. This yields the result as soon as the limits as the number of versions
grows and as $\eps\to0$ are suitably taken.
We define the two sets of versions now.

Let
$X^{(1)},X^{(2)},\ldots$ and $\tau^{(0)},\tau^{(1)},\ldots$ be
iid copies of $X$ and $\tau$ respectively. 
For $k\geq1$, let $\rk$ 
be defined as in~(\ref{eq:range}),
with $\xk$ replacing $X$.
Let now $\ZZ^{(1)}_n=\RR^{(1)}_n$ and for $k>1$  
\begin{equation}
\label{eq:zz}
\ZZ^{(k)}_n=\RR^{(k)}_n\setminus\left\{\cup_{i=1}^{k-1}\RR^{(i)}_n\right\}.
\end{equation}
We then define for each $N\geq1$
\begin{equation}
\label{eq:ttau} 
\tilde\tau^{(N)}=\{\tilde\tau_x^{(N)},\,x\in\Z^d\},
\end{equation}
where $\tilde\tau_x^{(N)}=\tau^{(k)}_x$, if $x\in\ZZ^{(k)}_N$ for some $k\geq1$, 
and $\tilde\tau_x^{(N)}=\tau^{(0)}_x$, otherwise. 

\begin{rmk}
 \label{rmk:tt}
$\tilde\tau^{(N)}$ and $\tau$ are equally distributed for every  $N\geq1$, whether or not $X^{(1)},X^{(2)},\ldots$ 
are given. 
In particular, $\tilde\tau^{(N)}$ is independent of $X^{(1)},X^{(2)},\ldots$.
\end{rmk}

Now let us consider two classes of clock processes $\ckn,\ck:\N\to[0,\infty)$, $k,N\geq1$:
\begin{equation}
\label{eq:cklock}
\ckn_n=\sum_{i=0}^n\tilde\tau^{(N)}_{\xk_i}\,T^{(k)}_i,\quad
\ck_n=\sum_{i=0}^n\tau^{(k)}_{\xk_i}\,T^{(k)}_i,\,\quad n\geq0,
\end{equation}
where $\{T^{(k)}_i,\,i\geq1\}=:T^{(k)},\,k\geq1$, are independent~families of iid mean 1 exponentials,
and their respective inverses $\ikn$ and $\ik$. Let then
\begin{equation}
\label{eq:yk}
\ykn_t=\tilde\tau^{(N)}_{\xk_{\ikn_t}},\quad \yk_t=\tau^{(k)}_{\xk_{\ik_t}},\quad t\geq0.
\end{equation}

\begin{rmk}
\label{rmk:ck} 
We remark that $\ckn=\ck=C$ and $\ykn=\yk=Y$ in distribution for all $k,N$; and that $\yk$, $k\geq1$, are iid; and 
$\ykn$, $k\geq1$, are iid given  $\tilde\tau^{(N)}$ for all $N$. 
The latter fact follows from the independence of
$\tilde\tau^{(N)}$ from $X^{(1)},X^{(2)},\ldots$ as remarked above (see Remark~\ref{rmk:tt}).
\end{rmk}

The reason why the argument we outlined above, and then started to fill the details of, works is that, as already remarked
and used in the previous section, the contributions to the processes come from only a few deep traps (in the sense
specified at the proof of Theorem~\ref{teo:conv1}, see~(\ref{eq:taudelta}) and~(\ref{eq:tauedelta}) above), which 
are far off one another, and, given Asssumption $C$, are unlikely to lie in any intersection of ranges. 
For this reason the difference between $\yk$ and $\ykn$ 
(suitably rescaled) come from sites which {\em are not} deep traps, and thus are negligible. 
In order to make this argument, let us fix $\eps,\d>0$ and define for each $k,N\geq1$
\begin{eqnarray}
\label{eq:tked}
\tked\=\{\tked_x:=\tk_x\1\{\tk_x>\d (\eps q_\eps)^{-1}\},\,x\in\Z^d\},\\
\label{eq:ttned}
\ttned\=\{\ttned_x:=\tilde\tau^{(N)}_x\1\{\tilde\tau^{(N)}_x>\d (\eps q_\eps)^{-1}\},\,x\in\Z^d\},
\end{eqnarray}
and let
\begin{equation}
\label{eq:cked}
\ckned_n=\sum_{i=0}^n\tilde\tau^{(N,\eps,\d)}_{\xk_i}\,T^{(k)}_i,\quad\cked_n=\sum_{i=0}^n\tau^{(k,\eps,\d)}_{\xk_i}\,T^{(k)}_i,\,\quad n\geq0,
\end{equation}
with $\ikned$ and $\iked$ their respective inverses, and
\begin{equation}\nn
\label{eq:yked}
\ykne_t=\eps q_\eps\tilde\tau^{(N)}_{\xk_{\ikn_{\eps^{-1}t}}},\quad \ykned_t=\eps q_\eps\tilde\tau^{(N)}_{\xk_{\ikned_{\eps^{-1}t}}},\,\quad 
\yke_t=\eps q_\eps\tau^{(k)}_{\xk_{\ik_{\eps^{-1}t}}},\quad \yked_t=\eps q_\eps\tau^{(k)}_{\xk_{\iked_{\eps^{-1}t}}}.
\end{equation}

We then have that $\ykne=\yke=\ye$ in distribution; $\yke$, $k\geq1$, are iid; and 
$\ykne$, $k\geq1$, are iid given   $\tilde\tau^{(N)}$ for all $N$.

Consider now
\begin{equation}
\label{eq:av1}
\frac1{K}\sum_{k=1}^{K}F(\ykne)=\frac1{K}\sum_{k=1}^{K}F(\ykned)
+\frac1{K}\sum_{k=1}^{K}\rkned,
\end{equation}
where $K$ is an arbitrary positive integer, and $\rkned=F(\ykne)-F(\ykned)$.

With an argument similar to the one giving~(\ref{eq:c6}), one finds that 
\begin{equation}
\label{eq:av2}
\lim_{\d\to0}\limsup_{\eps\to0}\rkned=0
\end{equation}
in probability for all $k$.


Since $\ckn=\ck=C$ in distribution, we have that, given $\eta>0$, there exists $S=S_K>0$ such that 
\begin{equation}
\label{eq:av3}
\P\left(\ckn_{\nu_\eps S}\leq\eps^{-1}T\right)=\P\left(\ck_{\nu_\eps S}\leq\eps^{-1}T\right)\leq\frac\eta{2K}
\end{equation}
for all $N\geq1$ and $\eps$ sufficiently small (see the second point of Remark~\ref{rmk:conv2} above).

From now on we take $N=N_\eps=2\rho_{\nu_\eps S}$.

For $k,\ell=1,\ldots,K$, let
\begin{equation}
\label{eq:av5}
A_{k,\ell}=\{\tk_{\xk_i}>\d (\eps q_\eps)^{-1}\mbox{ and }\xk_i\in\rl_{\nu_\eps S}\mbox{ for some }i\leq \nu_\eps S\},\quad
A_K=\cup_{k,\ell=1}^KA_{k,\ell}.
\end{equation}
Let also 
\begin{equation}
\label{eq:av6}
\I_{k,\ell}:=\I_{\nu_\eps S}(\xk,\xl)
\end{equation}
(see~(\ref{eq:int})).

Then, given $\lambda>0$
\begin{eqnarray}\nn
\P(A_{k,\ell})\=
\P\!\left(\sum_{x\in\I_{k,\ell}}\1\{\tk_x>\d (\eps q_\eps)^{-1}\}\geq1\right)\\
\label{eq:av7}
&\leq&\lambda N\P(\tau_0>\d (\eps q_\eps)^{-1})+\P(|\I_{k,\ell}|>\lambda N)+\P(R_{\nu_\eps S}^{(k)}>N).
\end{eqnarray}

By Assumption $B$, (\ref{eq:tail}, \ref{eq:ronr}, \ref{eq:vnu}),
the first term on the right of~(\ref{eq:av7}) is bounded above by
\begin{equation}
\label{eq:av8}
\lambda\, 3\d^{-\a} S^\a m\P(\tau_0>s_m)
\end{equation}
for all $\eps$ small enough, where $m=\rho_{\nu_\eps}$. By the definition of $s_m$ (see~(\ref{eq:sn}) above),
we may replace $m\P(\tau_0>s_m)$ by 1 in~(\ref{eq:av8}).

Putting this, Assumptions $A$ and~$C$ together, we conclude that $\P(A_{k,\ell})\to0$ as $\eps\to0$ for all $k,\ell$, and thus
\begin{equation}
\label{eq:av9}
\P(A_{K})\to0
\end{equation}
as $\eps\to0$ for every $K\geq1$.

We now go back to~(\ref{eq:av1}). By the above, the first term on its right can be written as
\begin{equation}
\label{eq:av10}
\frac1{K}\sum_{k=1}^{K}F(\yke)
+\frac1{K}\sum_{k=1}^{K}\rked
+\1\{A_K\cup B_K\}\frac1{K}\sum_{k=1}^{K}(F(\ykned)-F(\yked)),
\end{equation}
where
$\rked=F(\yked)-F(\yke)$ and
${\displaystyle B_K=\bigcup_{k=1}^K\left\{\{\ckn_{\nu_\eps S}\leq\eps^{-1}T\}\cup\{\ck_{\nu_\eps S}\leq\eps^{-1}T\}\right\}}$.
 This follows from the fact that outside $A_K\cup B_K$, we have that $\yked_t=\ykned_t$ for $t\in[0,T]$ and $1\leq k\leq K$, and 
all fixed $\eps,\delta,N$, as can be readily checked.

As in~(\ref{eq:av2}), we have that
\begin{equation}
\label{eq:av11}
\lim_{\d\to0}\limsup_{\eps\to0}\rked=0
\end{equation}
in probability for all $k$.

Now from (\ref{eq:av1},\ref{eq:av10}), and since  $F\in\B_u$, we get 
\begin{eqnarray}\nn
&\left|\esp\left[\left.F(\y1e)\right|\ttn\right]-\esp\left[F(Z)\right]\right|\leq&\\
&\left|\frac1{K}\sum_{k=1}^{K}\!\!\left(F(\ykne)-\esp\!\left[\left.F(\ykne)\right|\ttn\right]\right)\right|\!
+\!\left|\frac1{K}\sum_{k=1}^{K}\!\!\left(F(\yke)-\esp\!\left[F(\yke)\right]\right)\right|&
\label{eq:av12}
\end{eqnarray}
plus a term whose expectation is bounded above by 
\begin{equation}
\label{eq:av13}
\left|\esp(F(\y1e))-\esp(F(Z))\right|+\esp\left|\roned\right|+\esp\left|\r1ed\right|+2\|F\|_\infty\left(\P(A_K)+\P(B_K)\right),
\end{equation}
where we have used the fact that both $\esp\!\left[F(\yke)\right]$ and $\esp\!\left[\!\left.F(\ykne)\right|\ttn\right]$ are independent of $k$.

Recalling now Remark~\ref{rmk:ck}, and using Jensen, we find that the right hand side of~(\ref{eq:av12}) is bounded above by constant times
$K^{-1/2}$.
This and~(\ref{eq:s1}, \ref{eq:av2}, \ref{eq:av3}, \ref{eq:av9}, \ref{eq:av11}) yield
\begin{equation}
\label{eq:av14}
\limsup_{\eps\to0}\esp\left|\esp\left[\left.F(\y1e)\right|\ttn\right]-\esp\left[F(Z)\right]\right|\leq\mbox{const}\,(K^{-1/2}+\eta).
\end{equation}
Since $K$ and $\eta$ are arbitrary and the left hand side of~(\ref{eq:av14}) does not depend on either, we conclude that
\begin{equation}
\label{eq:av15}
\esp\left[\left.F(\y1e)\right|\ttn\right]\to\esp\left[F(Z)\right]
\end{equation}
in probability as $\eps\to0$, and the result follows from the fact that 
$\esp\left[\left.F(\y1e)\right|\ttn\right]$ and 
$\esp\left[\left.F(\ye)\right|\tau\right]$ have the same distribution for every $\eps>0,\, N\geq1$.
$\square$


\subsection{Stronger aging results}
\label{ssec:stra}

The above can be extended to strengthen the aging results of the previous section, under the same conditions of both this and that sections.

\begin{theo}
\label{teo:ages}
If $X$ is transient and Assumption $C$ holds, then we have that
\begin{eqnarray}
\label{eq:ages_r}
&\lim_{t\to\infty}\P(Y_t=Y_{t+\theta t}|\tau)=R(\theta),&\\
\label{eq:ages_pi}
&\lim_{t\to\infty}\P(Y_t=Y_{t+r}\,\mbox{ for all }r\in[0,\theta t/q_{1/t}]|\tau)=\Pi(\theta),&
\end{eqnarray}
in probability as $t\to\infty$, where $R$ and $\Pi$ are as in Theorem~\ref{teo:age} (and indeed, both equal the right hand side of~(\ref{eq:age3})
above).
\end{theo}

\noindent{\bf Sketch of proof}

An argument like that of Lemma~\ref{lm:conv1} can be used to get that
\begin{equation}
\label{eq:conve1a}
(\yoed_t)\to(\bzd_t)
\end{equation}
as $\eps\to0$ in distribution on $(D,J_1)$ (in here, differently from Lemma~\ref{lm:conv1} case, $\tau$ is integrated; the argument might of course
use 
a version of $\tau$ such that~(\ref{eq:conve1a}) holds for $\tau$ in a set of full probability). We may then extend Theorem~\ref{teo:str} with a
similar 
proof to get
\begin{equation}
\label{eq:str1}
\esp\left[\left.G(Y^{(1,\eps,\d)})\right|\tau\right]\to\esp\left[G(\bzd)\right],
\end{equation}
in probability as $\eps\to0$, for every $\d>0$, with $G:D\to\{0,1\}$ such that either
\begin{equation}
\label{eq:g}
G(U)=\1\{U_1=U_{1+r}\,\mbox{ for all }r\in[0,\theta]\}\mbox{ or }G(U)=\1\{U_1=U_{1+\theta}\},\,U\in D. 
\end{equation}
We further need to extend Lemma~\ref{lm:comp}
to get that for all $t>0$
\begin{equation}
\label{eq:str2}
\lim_{\delta\to0}\limsup_{\eps\to0}\P(Y^{(1,\eps)}_t\ne Y^{(1,\eps,\d)}_t)=0
\end{equation}
(here the proof can be made again by using special versions of $\tau$ like in the argument for Lemma~\ref{lm:comp}).

From~(\ref{eq:str1}) and~(\ref{eq:str2}), and since $\esp\left[G(\bzd)\right]\to\esp\left[G(Z)\right]$ as $\d\to0$, we get
\begin{equation}\label{eq:str3}
\lim_{\eps\to0}\esp\left[\left.G(Y^{(1,\eps)})\right|\tau\right]
=\lim_{\delta\to0}\lim_{\eps\to0}\esp\left[\left.G(Y^{(1,\eps,\d)})\right|\tau\right]
=\lim_{\delta\to0}\esp\left[\left.G(\bzd)\right|\tau\right]
=\esp\left[G(Z)\right],
\end{equation}
in probability as $\eps\to0$, for $G$ as in~(\ref{eq:g}).
$\square$

\begin{rmk}
\label{rmk:strqage}
Under the same conditions of Theorem~\ref{teo:ages}, and by the same reasoning as above, we have that
\begin{equation}
\label{eq:strqage}
\textstyle \lim_{t\to\infty}\P(\sup_{r\in[0,t]}Y_r<\sup_{r\in[0,t+\theta t]}Y_{r}|\tau)=\Omega(\theta)
\end{equation}
in probability. (See Remark~\ref{rmk:age4} above.)
\end{rmk}

\begin{rmk}
\label{rmk:striage}
With the extra condition of this section, argueing as in Subsection~\ref{ssec:iage} above, with the help of Theorem~\ref{teo:str}, 
we may establish stronger results for integrated aging functions. We state the following result as an example,
in the spirit of Subsection~\ref{ssec:iage} above.
\begin{theo}
\label{teo:striage}
If $X$ satisfies Assumptions  $A$, $B$ and $C$, then 
\begin{equation}
\label{eq:striage}
\lim_{\lambda,\mu\to\infty\atop{\lambda/\mu\to\theta}}\esp[\bar R(\lambda{\mathcal T},\mu {\mathcal T})|\tau]=R(\theta)
\end{equation}
in probability, with $R$ as in Theorem~\ref{teo:age} above.
\end{theo}

\end{rmk}

%% file: acks.tex
\vspace{.5cm}

\noindent{\bf Acknowledgements}

LRF would like to thank Marina Vachkovskaia for discussions at an early stage of this project. 
The authors thank an anonymous referee of an earlier version of this
paper for comments which led to substantial corrections and other improvements.